\documentclass[11pt]{article}
\pdfoutput=1
\usepackage{amsfonts,graphics}
 \usepackage{pgf,tikz,bbm,color}
\usepackage{epsfig}
\usepackage{times} 
\usepackage{url}
\usepackage{float}
\newfloat{algorithm}{tbp}{lop}
\floatname{algorithm}{Algorithm}

\def\sq{ { { \hfill\rule{3mm}{3mm} } } }

\newfont{\Bb}{msbm10 scaled\magstep0}
\def\IR{\mbox {\Bb R}}

\def\IN{\mbox {\Bb N}}

\newcommand{\comm}[1]{}

\newcommand{\bA}{{\bf A}}

\newcommand{\bC}{{\bf C}}

\newcommand{\bE}{{\bf E}}
\newcommand{\bF}{{\bf F}}

\newcommand{\bP}{{\bf P}}
\newcommand{\bS}{{\bf S}}
\newcommand{\bY}{{\bf Y}}

\newcommand{\bL}{{\bf L}}
\newcommand{\bM}{{\bf M}}

\newcommand{\bI}{{\bf I}}

\newcommand{\bQ}{{\bf Q}}
\newcommand{\bW}{{\bf W}}
\newcommand{\bR}{{\bf R}}
\newcommand{\bT}{{\bf T}}
\newcommand{\bZ}{{\bf Z}}

\newcommand{\bV}{{\bf V}}
\newcommand{\bU}{{\bf U}}
\newcommand{\bfa}{{\bf a}}
\newcommand{\bb}{{\bf b}}
\newcommand{\bc}{{\bf c}}

\newcommand{\bfe}{{\bf e}}
\newcommand{\bff}{{\bf f}}
\newcommand{\bp}{{\bf p}}
\newcommand{\bq}{{\bf q}}
\newcommand{\br}{{\bf r}}

\newcommand{\bv}{{\bf v}}
\newcommand{\bw}{{\bf w}}
\newcommand{\bx}{{\bf x}}
\newcommand{\by}{{\bf y}}
\newcommand{\bz}{{\bf z}}

\newcommand{\bPhi}{\mbox{\boldmath$\Phi$}}
\newcommand{\bPsi}{\mbox{\boldmath$\Psi$}}
\newcommand{\bOmega}{\mbox{\boldmath$\Omega$}}
\def\bell{\mbox{\boldmath$\ell$}}

\newcommand{\cO}{ {\cal O} }

\newcommand{\cP}{ {\cal P} }

\newcommand{\cQ}{ {\cal Q} }

\newcommand{\tol} {\mbox{{\it tol}}}
\mathcode`\@="8000
{\catcode`\@=\active \gdef@{\mkern1mu}} 
\newcommand{\Ran}{{\rm Ran}}

\newcommand{\be}{\begin{equation}}
\newcommand{\ee}{\end{equation}}
\newcommand{\ba}{\begin{array}}
\newcommand{\ea}{\end{array}}

\def\l1{\ell_1}
\def\l2{\ell_2}

\def\L1{{\cal L}_1}
\def\L2{{\cal L}_2}

\def\h1{{\it h}_1}
\def\h2{{\it h}_2}

\def\H1{{\cal H}_1}
\def\H2{{\cal H}_2}

\def\H2{ { {\cal H}_2} }

\def\sq{ { { \hfill\rule{3mm}{3mm} } } }
\def\rank{ { {{\rm rank }\,} } }

\def\QQ{\bQ}
\def\qq{\bq}
\def\RR{\bR}
\def\rr{\br}

\newtheorem{theorem}{Theorem}[section]
\newtheorem{lemma}{Lemma}[section]
\newtheorem{definition}{Definition}[section]

\newcommand{\blem}{\begin{lemma}}
\newcommand{\elem}{\end{lemma}}
\newcommand{\bea}{\begin{eqnarray*}}
\newcommand{\eea}{\end{eqnarray*}}
\newcommand{\bean}{\begin{eqnarray}}
\newcommand{\eean}{\end{eqnarray}}

\makeatletter
\@addtoreset{equation}{section}
\makeatother

\title{\Large{\bf A DEIM Induced CUR Factorization
}\thanks{
This work was supported in part by
AFOSR grant FA9550-12-1-0155 and by NSF grant CCF-1320866.
}}

\author{{\sc D. C. Sorensen}\footnote{Department of Computational 
and Applied Mathematics, MS 134, Rice University, Houston, Texas 77005-1892.}
{\sc\ and M. Embree}\footnote{Department of Mathematics, 225 Stanger Street 0123, Virginia Tech, Blacksburg, Virginia 24061} 
\\[2ex]
{\tt e-mail: sorensen@rice.edu, embree@vt.edu }
\vspace{.2cm}
}
\date{\today}

\begin{document}
\maketitle

\begin{abstract}
We derive a CUR approximate matrix factorization
based on the Discrete Empirical Interpolation Method (DEIM).\ \ 
For a given matrix $\bA$, such a factorization provides a low rank approximate decomposition of the form $\bA \approx \bC \bU \bR$,
where $\bC$ and $\bR$ are subsets of the columns and rows of $\bA$, and $\bU$ is constructed to make  $\bC\bU \bR $  a good approximation. 
Given a low-rank singular value decomposition $\bA \approx \bV \bS \bW^T$, the DEIM procedure uses $\bV$ and $\bW$ to select the columns and rows of $\bA$ that form $\bC$ and $\bR$.\ \ 
Through an error analysis applicable to a general class of CUR factorizations, we show that the accuracy tracks the optimal approximation error within a factor that depends on the conditioning of submatrices of $\bV$ and $\bW$.\ \ For very large problems, $\bV$ and $\bW$ can be approximated well using an incremental QR algorithm that makes only one pass through $\bA$.\ \  Numerical examples illustrate the favorable performance of the DEIM-CUR method compared to CUR approximations based on leverage scores.
\end{abstract}

\section{Introduction}
This work presents a new CUR matrix factorization
based upon the Discrete Empirical Interpolation Method (DEIM).\ \ 
A CUR factorization is a low rank approximation
of a matrix $\bA \in \IR^{m \times n}$ of the form $\bA \approx \bC \bU \bR$, where
$\bC = \bA(:,\bq) \in \IR^{m \times k}$ is a subset of the columns of $\bA$ and
$\bR = \bA(\bp,:) \in \IR^{k \times n} $ is a subset of the rows of $\bA$.\ \  
(We generally assume $m\ge n$ throughout.)
The $k \times k$ matrix $\bU$ 
is constructed to assure that $\bC \bU \bR$ is a good approximation to  $\bA $.
Assuming the best rank-$k$ singular value decomposition (SVD) $\bA \approx \bV \bS \bW^T$ is available, the algorithm uses the DEIM index selection procedure, $\bq = {\rm DEIM}(\bV)$ and $\bp = {\rm DEIM}(\bW)$, to determine $\bC$ and $\bR$.
The resulting approximate factorization is nearly as accurate as the best rank-$k$ SVD, with
\[
  \| \bA - \bC \bU \bR \| \le (\eta_p + \eta_q)\,\sigma_{k+1},
\]
where $\sigma_{k+1}$ is the first neglected singular value of $\bA$,  $\eta_p \equiv \| \bV(\bp,:\,)^{-1} \|$, and $\eta_q \equiv \| \bW(\bq,:\,)^{-1} \|$.  

Here and throughout, $\|\cdot\|$ denotes the vector 2-norm and the matrix norm it induces, and $\|\cdot\|_F$ is the Frobenius norm.   We use MATLAB notation to index vectors and matrices, so that, e.g., $\bA(\bp,:)$ denotes the $k$ rows of $\bA$ whose indices are specified by the entries of the vector $\bp\in\IN^k$, while $\bA(:,\bq)$ denotes the $k$ columns of $\bA$ indexed by $\bq\in\IN^k$.

\begin{figure}[b!]
\begin{center}
\includegraphics[scale=0.475]{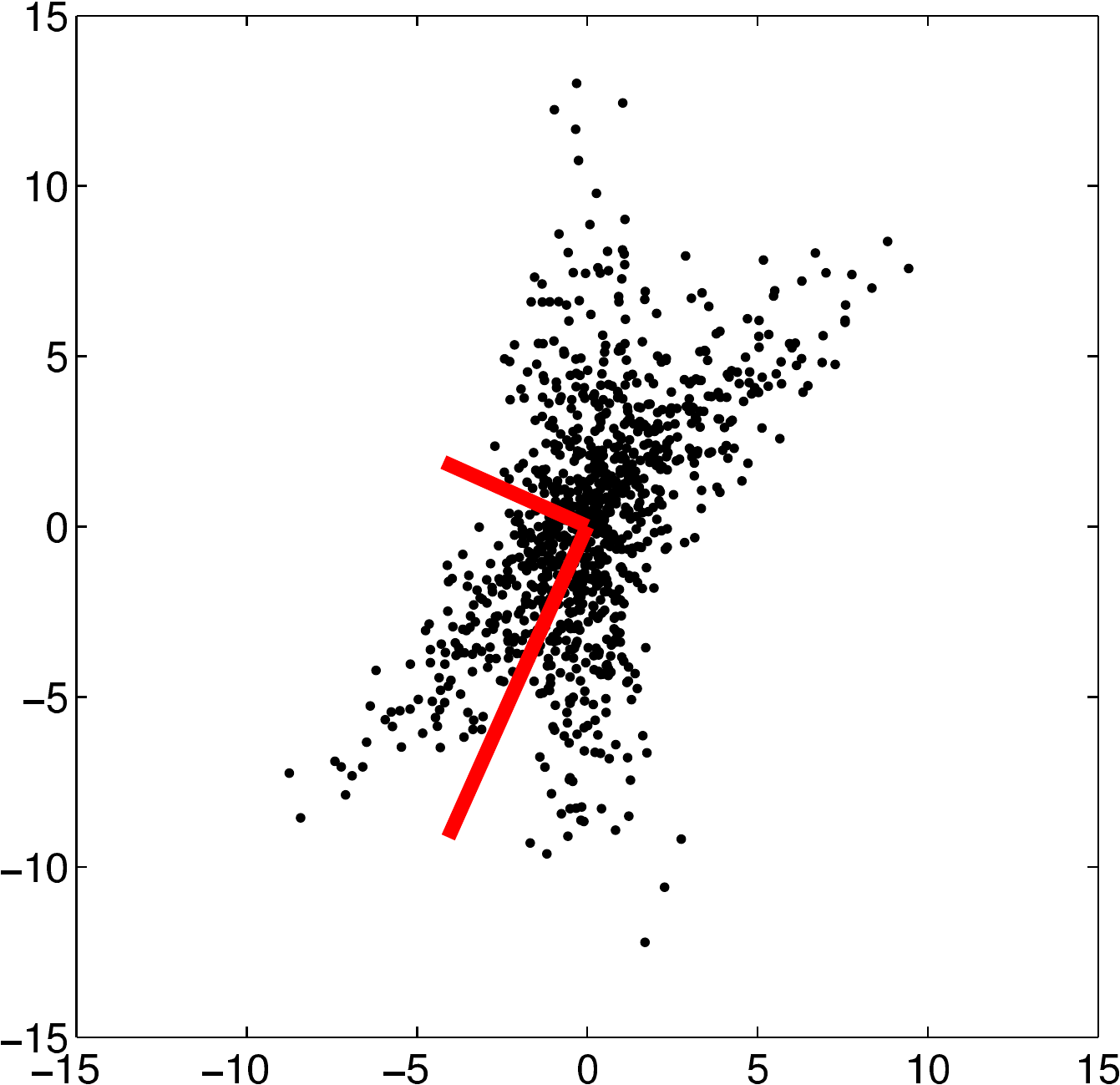}\qquad\quad
\includegraphics[scale=0.475]{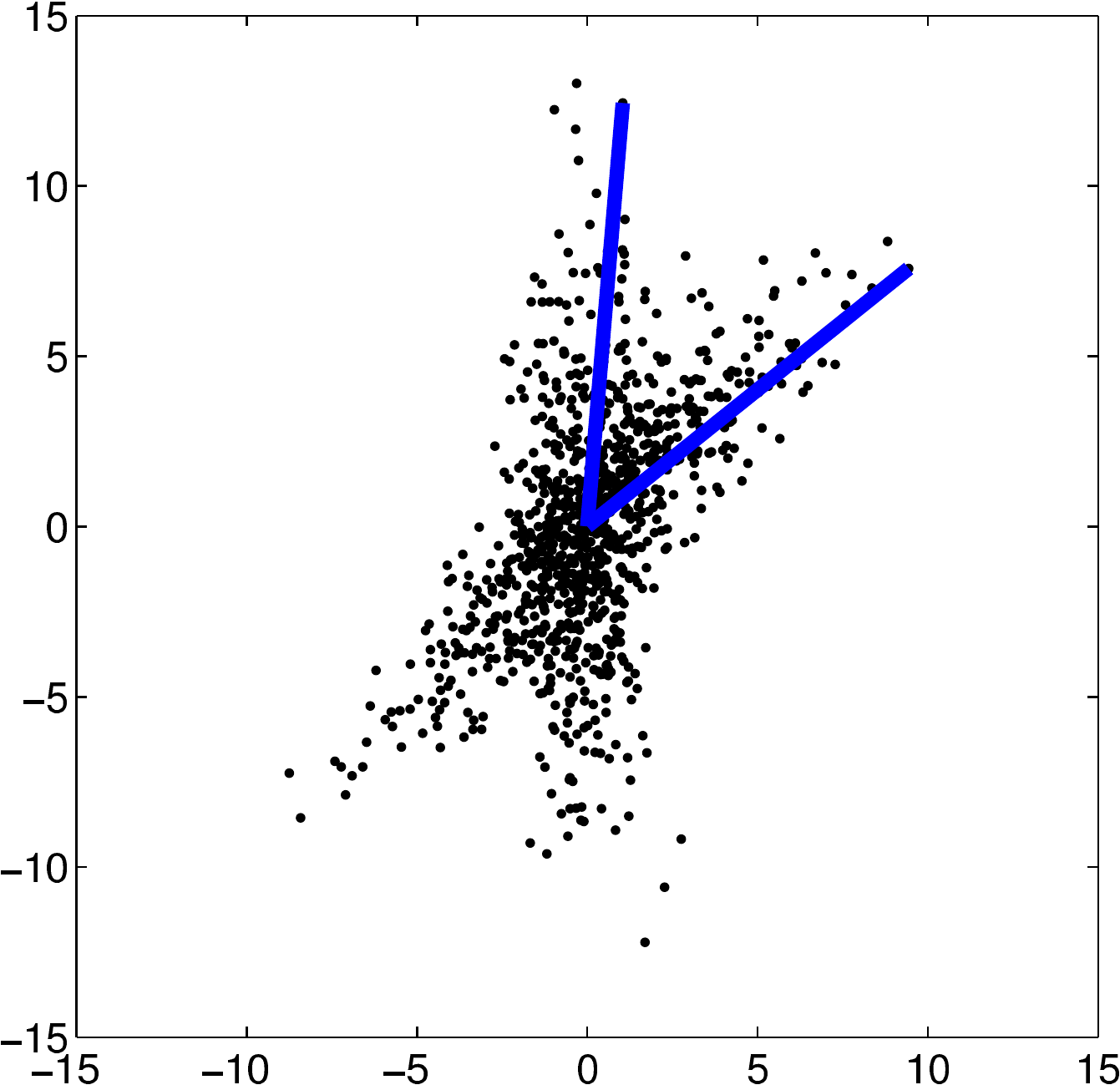}
\end{center}

\vspace*{-10pt}
\caption{\label{fig:mv_example} Comparison of singular vectors (left, scaled, in red) 
and DEIM-CUR columns (right, in blue) for a data set drawn from two multivariate normal distributions having different principal axes.}
\end{figure}

The CUR factorization is an important tool for handling large-scale data sets, offering two advantages over the SVD:  when $\bA$ is sparse, so too are $\bC$ and $\bR$, unlike the matrices $\bV$ and $\bW$ of singular vectors; and the columns and rows that comprise $\bC$ and $\bR$ are representative of the data (e.g., sparse, nonnegative, integer valued, etc.).\ \  The following simple example, adapted from Mahoney and Drineas~\cite[Fig.~1b]{MD09}, illustrates the latter advantage.  Construct $\bA\in\IR^{2\times n}$ so that its first $n/2$ columns have the form
\[ \left[\begin{array}{c} x_1 \\ x_2 \end{array}\right]\]
and the remaining $n/2$ columns have the form
\[    {\sqrt{2}\over2} \left[\begin{array}{cc} -1 & 1 \\ 1 & 1\end{array}\right] 
    \left[\begin{array}{c} x_1 \\ x_2 \end{array}\right],
\]
where in both cases $x_1 \sim N(0,1)$ and $x_2 \sim N(0,4^2)$ are independent samples of normal random variables, i.e.,  
the columns of $\bA$ are drawn from two different multivariate normal distributions.  Figure~\ref{fig:mv_example} shows that the two left singular vectors, though orthogonal by construction, fail to represent the true nature of the data; in contrast, the first two columns selected by the DEIM-CUR procedure give a much better overall representation.  While trivial in this two-dimensional case, one can imagine the utility of such approximations for high-dimensional data.  We shall illustrate the advantages of CUR approximations with further computational examples in Section~\ref{sec:examples}.

CUR-type factorizations originated with ``pseudoskeleton'' approximations~\cite{GTZ97} and pivoted, truncated QR decompositions~\cite{Ste99}; in recent years many new algorithms have been proposed in the numerical linear algebra and theoretical computer science literatures. 
Some approaches seek to maximize the \emph{volume} of the decomposition~\cite{GTZ97,TKB12}.  
Numerous other algorithms instead use \emph{leverage scores}~\cite{BW,DMM08,MD09,WZ13}.  These methods typically first compute a singular value decomposition%
\footnote{We use the nonstandard notation $\bV \bS \bW^T$ for the SVD to avoid conflicts with $\bU$ in the standard CUR notation.} 
$\bA = \bV\bS\bW^T$ (or an approximation to it), with $\bV\in\IR^{m\times n}$, $\bW\in\IR^{n\times n}$.\ \ 
The leverage score for the $j$th row ($k$th column) of $\bA$ is the squared two-norm of the $j$th row of $\bV$ ($k$th row of $\bW$).  When scaled by the number of singular vectors, these leverage scores give probability distributions for randomly sampling the columns and rows to form $\bC$ and $\bR$.  This approach leads to probabilistic bounds on $\|\bA - \bC \bU \bR \|_F$~\cite{DMM08,MD09}.  In cases where $\bA$ has small singular values (precisely the case where one would seek a low-rank factorization), the singular vectors can be sensitive to perturbations to $\bA$, making the leverages scores unstable~\cite{IW}.  Thus leverage scores are often computed using only the leading few singular vectors, but the choice of how many vectors to keep can be somewhat ad hoc.

The algorithm described in Sections~\ref{sec:factor} and~\ref{sec:deim} is entirely deterministic and involves few (if any) parameters.  The method is supported by an error analysis in Section~\ref{sec:theory} that also applies to a broad class of CUR factorizations.  
This section includes an improved bound on the error constants 
$\eta_p$ and $\eta_q$ for DEIM row and column selection, which also applies to  the analysis of DEIM-based model order reduction~\cite{CS10}.
In Section~\ref{sec:qr} we propose a novel incremental QR algorithm for approximating the SVD
(and potentially also approximating leverage scores).  Section~\ref{sec:examples} illustrates the performance of this new CUR factorization on several examples.

In many applications one cares primarily about key columns or rows of $\bA$, rather than an explicit $\bA=\bC\bU\bR$ factorization.
The DEIM technique, which identifies rows and columns of $\bA$ independently, can easily be used to 
select only columns or rows, leading to an ``interpolatory decomposition'' of the form 
$\bA = \bC\widehat{\bU}$ or $\bA = \widehat{\bU}\bR$;
such factorizations have the advantage that $\widehat{\bU}$ can be much better conditioned than the $\bU$ matrix in the CUR factorization.  For further details about general interpolatory decompositions, see~\cite[\S1]{CGMR05}.

\section{CUR Factorization} \label{sec:factor}
We are concerned with large matrices $\bA\in\IR^{m\times n}$ that represent nearly low-rank data, which can therefore be expressed as
\begin{equation} \label{eq:CURF}
  \bA = \bC \bU \bR  + \bF,
\end{equation}
with $\|\bF\|$ small relative to $\|\bA\|$.\ \ 
The matrix $\bC\in\IR^{m\times k}$ is formed by extracting $k$ columns from $\bA$, and $\bR\in\IR^{k\times n}$ from $k$ rows of $\bA$.  The selected row and column indices are
stored in the vectors $\bp, \bq\in\IN^k$, so that $\bC = \bA(:,\bq)$ and $\bR = \bA(\bp,:)$.  
Our choice for $\bp$ and $\bq$ is guided by knowledge of the  rank-$k$ SVD (or an approximation to it).  Before detailing the method for selecting these indices, we discuss how, given $\bp$ and $\bq$, one should construct $\bU$ so that $\bC\bU\bR$ satisfies desirable approximation properties.

As motivation, suppose for the moment that $\bA$ has exact rank~$k$, and $\bC$ and $\bR$ are full-rank subsets of the columns and rows of $\bA$.\ \ 
Now let $\bY\in\IR^{m\times k}$ and $\bZ\in\IR^{n\times k}$ be any matrices that satisfy $\bY^T \bC = \bR \bZ = \bI\in\IR^{k\times k}$.  Then $\bC\bY^T$ is a projector onto $\Ran(\bC)=\Ran(\bA)$ and $(\bZ \bR)^T$ is a projector onto $\Ran(\bR^T)=\Ran(\bA^T)$, where $\Ran(\cdot)$ denotes the range (column space).  It follows that $\bC\bY^T \bA = \bA$ and $(\bZ\bR)^T\bA^T = \bA^T$.\ \   Putting $\bU \equiv \bY^T\bA\bZ$ gives 
\[ \bC\bU\bR = \bC\bY^T\bA\bZ\bR = \bA\bZ\bR = \bA.\]
Thus, any choice of $\bY$ and $\bZ$ that satisfies $\bY^T \bC = \bR \bZ = \bI$ gives a $\bU$ such that  $\bC\bU\bR$ exactly recovers $\bA$.\ \ 
In general different choices for $\bY$ and $\bZ$ give different $\bU = \bY^T \bA \bZ$.

Now consider the general case~(\ref{eq:CURF}).  Once $\bp$, $\bq$, $\bY$, and $\bZ$ have been specified, then
\[
   \bU = \bY^T \bA \bZ  \quad {\rm and } \quad  \bF \equiv \bA - \bC \bU \bR.
\]   
One might design $\bY$ and $\bZ$ so that $\bC\bU\bR$ matches the selected columns $\bC = \bA(:,\bq)$ 
and rows $\bR = \bA(\bp,:)$ of $\bA$ exactly.   This can be accomplished with \emph{interpolatory projectors}, which we discuss in detail in the next section.  For now, let $\bP = \bI(@:@,\bp)\in\IR^{m\times k}$ and $\bQ = \bI(@:@,\bq)\in\IR^{n\times k}$ be submatrices of the identity, so that $\bP^T \bfa = \bfa(\bp)$ and $\bb^T \bQ = \bb(\bq)^T$ for arbitrary vectors
$\bfa$ and $\bb$ of appropriate dimensions. 
Now define $\bY^T = (\bP^T \bC)^{-1} \bP^T$ and $\bZ = \bQ (\bR\bQ)^{-1}$ (presuming $\bP^T\bC$ and $\bR\bQ$ are invertible).  Then since $\bC = \bA(:,\bq)$ and $\bR = \bA(\bp,:)$,
\[
   \bP^T \bC = \bC(\bp,:) = \bA(\bp,\bq) 
   \quad \mbox{and}\quad
   \bR \bQ = \bR(:,\bq) = \bA(\bp,\bq),\]
so
\[ \bU = \bY^T \bA \bZ 
       = (\bP^T\bC)^{-1}\bP^T \bA\bQ (\bR\bQ)^{-1} = \bA(\bp,\bq)^{-1}\bA(\bp,\bq) \bA(\bp,\bq)^{-1} = \bA(\bp,\bq)^{-1}.\]
This CUR approximation matches the $\bq$ columns and $\bp$ rows of $\bA$,
\[
   \bA(:,\bq) = \bC \bU \bR(:,\bq) \ \ {\rm and} \ \
   \bA(\bp,:) = \bC(\bp,:) \bU \bR,
\]
and, in our experiments, usually delivers a very good approximation.
However, a CUR factorization with better theoretical approximation properties
results from orthogonal projection, as originally suggested by Stewart~\cite[p.~320]{Ste99}; see also, e.g., Mahoney and Drineas \cite{MD09}.
Given a selection of indices $\bp$ and $\bq$, again put
\[
   \bC = \bA(:,\bq) \quad\mbox{and}\quad \bR = \bA(\bp,:).
\]
Assume that $\bC$ and $\bR$ 
 both have full rank $k$, and now let 
$\bY^T = \bC^I \equiv (\bC^T \bC)^{-1} \bC^T $ 
and $\bZ = \bR^I \equiv \bR^T (\bR \bR^T)^{-1} $ denote left and right inverses of $\bC$ and $\bR$.  
These choices also satisfy $\bY^T\bC = \bI$ and $\bR\bZ=\bI$, but now $\bC\bY^T = \bC\bC^I$ and 
$\bZ\bR = \bR^I\bR$ are \emph{orthogonal projectors}. 
We compute
\[
   \bU = \bY^T \bA\bZ = \bC^I \bA \bR^I,
\]
yielding a CUR factorization that can be viewed as a two step process: 
first the columns of $\bA$ are projected onto Ran($\bC$), 
then the result is projected onto the row space of $\bR$:
\[
   1) \  \  \bM = \bC\bC^I \bA, \quad 2)\ \  \bC \bU \bR  =  \bM \bR^I \bR.
\] 
Both steps are optimal with respect to the 2-norm error, which is the primary
source of the excellent approximation properties of this approach.

Several strategies for selecting $\bp$ and $\bq$ have been
proposed.%
\footnote{In the theoretical computer science literature, one often takes $\bC$ and/or $\bR$ to have rank larger than $k$, but then builds $\bU$ with rank $k$.  By selecting these extra columns and/or rows, one seeks to get within some factor $1+\varepsilon$ of the optimal approximation; see, e.g.,~\cite{BW}.}
The approach presented in the next section is simple to implement
and has complexity $ m@k $ and $n@k $ to select the
indices $\bp$ and $\bq$, provided the leading $k$ right and left singular vectors of $\bA$ are available.
Thus the overall complexity is dominated by the construction of the
rank-$k$ SVD  $\bA \approx \bV \bS \bW^T$, where $\bV^T \bV = \bW^T \bW = \bI\in\IR^{k\times k}$
and $\bS = {\rm diag}(\sigma_1, \sigma_2, \ldots , \sigma_{k})$ is the
$k \times k$ matrix of dominant singular values $\sigma_1 \ge \sigma_2  \ge \cdots \ge \sigma_{k}$.

\section{DEIM} \label{sec:deim}
The DEIM point selection algorithm was first presented in~\cite{CS10} in the
context of model order reduction for nonlinear dynamical systems, and is a discrete variant of the Empirical Interpolation Method originally proposed in~\cite{BMNP}.   
The DEIM procedure operates on the singular vector matrices $\bV$ and $\bW$ independently to select the row indices $\bp$ and column indices $\bq$.  We explain the process for selecting $\bp$; applying the same steps to $\bW$ yields~$\bq$.
To derive the method, we elaborate upon the interpolatory projectors introduced in the last section.

\begin{definition}
Given a full rank matrix $\bV \in \IR^{m\times k}$ and a set of distinct indices $\bp\in\IN^k$, the \emph{interpolatory projector} for $\bp$ onto $\Ran(\bV)$ is
\begin{equation} \label{eq:intproj}
\cP \equiv \bV(\bP^T\bV)^{-1}\bP^T,
\end{equation}
where $\bP = \bI(@:@,\bp)\in\IR^{m\times k}$, 
provided $\bP^T\bV$ is invertible.
\end{definition}

In general $\cP$ is an oblique projector, and it has an important property not generally enjoyed by orthogonal projectors:  for any  $\bx\in\IR^m$, 
\[ (\cP\bx)(\bp) \,=\, \bP^T\cP\bx \,=\, \bP^T\bV(\bP^T\bV)^{-1}\bP^T\bx \,=\, \bP^T\bx \,=\, \bx(\bp),\]
so the projected vector $\cP\bx$ matches $\bx$ in the $\bp$ entries, justifying the name ``interpolatory projector.''

The DEIM algorithm processes the columns of 
\[ \bV = \left[\begin{array}{cccc}\bv_1 & \bv_2 & \cdots & \bv_k\end{array}\right] \]
one at a time, starting from the leading singular vector $\bv_1$.  Each step processes the next singular vector to produce the next index.  The first index $p_1$ corresponds to the largest magnitude entry in $\bv_1$:
\[ |\bv_1(p_1)| = \|\bv_1\|_\infty.\]
Now define $\bp_1 \equiv [p_1]$, and let
\[ \cP_1 \equiv \bv_1 (\bP_1^T\bv_1)^{-1}\bP_1^T\]
denote the interpolatory projector for $\bp_1$ onto $\Ran(\bv_1)$.
The second index $p_2$ corresponds to the largest entry in $\bv_2$, after the interpolatory projection in the $\bv_1$ direction has been removed:
\begin{eqnarray*}
 \br_2\!\! &\equiv&\!\! \bv_2 - \cP_1 \bv_2\\[.25em]
 |\br_2(p_2)|\!\! &=&\!\! \|\br_2\|_\infty.
\end{eqnarray*}
Notice that $\br_2(p_1) = 0$, since $\cP_1\bv_2$ matches $\bv_2$ in the $p_1$ position, a consequence of interpolatory projection.  This property ensures the process will never produce duplicate indices.

Now suppose we have $j-1$ indices, with
\[ \bp_{j-1} \equiv \left[\begin{array}{c} p_1 \\ \vdots \\ p_{j-1}\end{array}\right], \quad
    \bP_{j-1} \equiv \bI(@:@,\bp_{j-1}), \quad
     \bV_{j-1} \equiv [\begin{array}{ccc} \bv_1 & \cdots & \bv_{j-1}\end{array}], \quad
    \cP_{j-1} \equiv \bV_{j-1}^{}(\bP_{j-1}^T\bV_{j-1}^{})^{-1}\bP_{j-1}^T.
    \]
To select $p_{j}$, remove from $\bv_j$ its interpolatory projection onto indices $\bp_{j-1}$ and take the largest remaining entry:
\begin{eqnarray*}
 \br_{j} \!\! &\equiv&\!\! \bv_j - \cP_{j-1} \bv_j\\[.25em]
 |\br_{j}(p_{j})|\!\! &=&\!\! \|\br_{j}\|_\infty.
\end{eqnarray*}
Implementations should not explicitly construct these projectors; see the pseudocode in Algorithm~\ref{fig:Deim_Arnoldi} for details.

Those familiar with partially pivoted LU decomposition will notice,
on a moment's reflection, 
that this index selection scheme 
is exactly equivalent to the index selection of partial pivoting.   
This arrangement is equivalent to the ``left looking" variant of
LU factorization~\cite[sect.~5.4]{DDSV98}, but with two important differences.  
First, there are no explicit row interchanges in DEIM, as there are
in LU factorization.   Second, the original basis vectors 
(columns of $\bV$) are not replaced
with the residual vectors, as happens in traditional LU decomposition.  
(In the context of model reduction,
it is preferable to keep the nice orthogonal basis intact for use 
as a reduced basis.)
We will exploit this connection with partially pivoted LU factorization 
to analyze the approximation properties of DEIM.

Since the DEIM algorithm processes the singular
vectors sequentially, from most to least significant, it 
 introduces new singular vector information in a 
coherent manner as it successively selects the $k$~indices.  Contrast 
this to index selection strategies based on 
leverage scores, where all singular vectors are 
incorporated at once via row norms of $\bV$ and $\bW$; to account for the fact that higher singular vectors are less significant, such approaches often instead compute leverage scores using only a few of the leading singular vectors.%
\footnote{A potential limitation of the DEIM approach is that $\br_j$ could have multiple entries that have nearly the same magnitude, but only one index is selected at the $j$th step; if the other large-magnitude entries in $\br_j$ are not significant in subsequent $\br_\ell$ vectors, the corresponding indices will not be selected.  One can imagine modifications of the selection algorithm to account for such situations, e.g., by processing multiple singular vectors at a time.}

\begin{algorithm}[t!]
\begin{center}
\fbox{
\begin{tabular}{p{10cm}}
\ \\[-25pt]
\begin{tabbing}
defghasijkl\=bbb\=ccc\= \kill
$\ \ \ $ {\bf Input:} \> $ \bV$,  an $m \times k$ matrix ($m\ge k$)\\[5pt]
$\  $  {\bf Output:}  \>  $\bp$, an integer vector with $k$ distinct entries in $ \{1,\ldots,m \}$ \\[7pt]
$\ \ \ $          \> $ \bv = \bV(:,1)$  \\
$\ \ \ $          \> $[\sim , p_1] = {\rm max}( | \bv |)$\\
$\ \ \ $           \>$ \bp = [p_1] $ \\
$\ \ \ $           \>{\bf for} $j=2,3,\ldots, k$  \\
$\ \ \ $               \>\> $ \bv = \bV(:,j)$  \\
$\ \ \ $               \>\> $ \bc =  \bV(\bp,1:j-1)^{-1} \bv(\bp)$  \\
$\ \ \ $               \>\> $ \br = \bv - \bV(:,1:j-1) \bc$  \\
$\ \ \ $               \>\>$[\sim , p_j] = {\rm max}(| \br |)$\\ 
$\ \ \ $               \>\>$ \bp = [\bp;\ p_j] $ \\
$\ \ \ $             \>{\bf end}\\
\end{tabbing}\\[-20pt]\end{tabular}}
\end{center}

\vspace*{-7pt}
\caption {DEIM point selection algorithm.}
\label{fig:Deim_Arnoldi}
\end{algorithm}


For the interpolatory projector $\cP_j$ to exist at the $j$th step, $\bP_{j-1}^T\bV_{j-1}$ must be nonsingular.  The linear independence of the columns of $\bV$ assures this.  In the following, $\bfe_j$ denotes the $j$th column of the identity matrix.

\begin{lemma} \label{DEIM_inv}
Let $\bP_j = [ \bfe_{p_1}, \bfe_{p_2}, \ldots, \bfe_{p_j}]$ and
let $\bV_j = [ \bv_1, \bv_2, \ldots, \bv_j]$ for $1 \le j \le k$.
If $\rank(\bV) = k$, then $\bP_j^T \bV_j$ is nonsingular for $1 \le j \le k$.
\end{lemma}
\noindent{\it Proof:} 
Suppose $\bP_{j-1}^T \bV_{j-1} $ is nonsingular and
let $\br_j = \bv_j - \bV_{j-1}(\bP_{j-1}^T \bV_{j-1})^{-1} \bP_{j-1}^T \bv_j.$
Then $\|\br_j\|_{\infty}  > 0$, for otherwise
$ {\bf 0} = \bv_j - \bV_{j-1} \bc_{j-1}$, in violation of the assumption that $\rank (\bV) = k$.
Thus
\begin{equation}
\label{eta_formula}
   0 < |\bfe_{p_j}^T \br_j| =| \bfe_{p_j}^T \bv_j - \bfe_{p_j}^T \bV_{j-1}(\bP_{j-1}^T \bV_{j-1})^{-1} \bP_{j-1}^T \bv_j |,
\end{equation}
where $p_j$ is the $j$th DEIM interpolation point.
Now factor
\begin{equation}
\label{PV_formula}
   \bP_j^T \bV_j =  
                   \left[ \begin{array}{cc} \bP_{j-1}^T \bV_{j-1} & \bP_{j-1}^T\bv_j \\
                                             \bfe_{p_j}^T \bV_{j-1} & \bfe_{p_j}^T \bv_j  
                            \end{array}
                     \right]=
                   \left[ \begin{array}{cc} \bI_{j-1} & {\bf 0} \\
                                             \bfe_{p_j}^T \bV_{j-1}(\bP_{j-1}^T \bV_{j-1})^{-1} & 1  
                            \end{array}
                     \right]
                   \left[ \begin{array}{cc} \bP_{j-1}^T \bV_{j-1} & \bP_{j-1}^T\bv_j \\
                                             {\bf 0} & \nu_j  
                            \end{array}
                     \right],
\end{equation}
where 
\[
   \nu_j = \bfe_{p_j}^T \bv_j - \bfe_{p_j}^T \bV_{j-1}(\bP_{j-1}^T \bV_{j-1})^{-1} \bP_{j-1}^T \bv_j .
\]
The inequality~(\ref{eta_formula}) implies $\nu_j \ne 0$ and hence equation~(\ref{PV_formula})
implies $\bP_j^T \bV_j$ is nonsingular.  Since $\bfe_{p_1}^T\bv_1 \ne 0$, this argument
provides an inductive proof that  $\bP_j^T \bV_j $  is nonsingular for  $1 \le j \le k$.  \sq

\section{CUR Approximation Properties} \label{sec:theory}

While the theory presented in this section was designed to
bound $\|\bA-\bC\bU\bR\|$ for the DEIM-CUR method, 
the analysis applies to \emph{any} CUR factorization with 
full rank $\bC \in \IR^{m\times k}$ and 
$\bR\in \IR^{k\times n}$, and $\bU = \bC^I\bA\bR^I$, 
regardless of the procedure used for selecting the columns 
and rows.%
\footnote{We are grateful to Ilse Ipsen for noting
the applicability of this analysis to all such $\bC\bU\bR$ factorizations,
and for also pointing out that, given knowledge of \emph{all} the singular values and vectors of $\bA$, our Lemma~\ref{ProjErr} can be sharpened via application of~\cite[Thm.~9.1]{HMT11}.  Indeed, Ipsen observes that the interpolatory projector proof of Lemma~\ref{ProjErr} can be adapted to simplify the multipage proof of \cite[Thm.~9.1]{HMT11}.}

Consider a CUR factorization that uses row indices $\bp \in \IN^k$ and column indices $\bq\in\IN^k$, and set
\[ \bP = \bI(@@:@@,\bp) = [\bfe_{p_1}, \ldots, \bfe_{p_k}]\in\IR^{m\times k}, \qquad
    \bQ = \bI(@@:@@,\bq) = [\bfe_{q_1}, \ldots, \bfe_{q_k}]\in\IR^{n\times k}.
 \]
The first step in this analysis bounds the mismatch
between $\bA$ and its interpolatory projection $\cP\bA$.


\begin{lemma}
\label{DEIM_bound}
Assume $\bP^T\bV$ is invertible and let $\cP = \bV (\bP^T\bV)^{-1}\bP^T$ be the interpolatory
projector~$(\ref{eq:intproj})$. If $\bV^T\bV=\bI$,  then  any $\bA\in\IR^{m\times n}$ satisfies
\[
  \| \bA - \cP\bA \| \le \|(\bP^T \bV)^{-1} \| \|(\bI - \bV \bV^T)\bA \|.
\]
Additionally, if $\bV$ consists of the leading $k$ left singular vectors of $\bA$, then
\[
  \| \bA - \cP\bA \|  =  \|(\bI - \cP) \bA\| \le \|(\bP^T \bV)^{-1} \|\, \sigma_{k+1}.
\]
\end{lemma}

\noindent{\it Proof:} First note 
that $\cP \bV = \bV (\bP^T \bV)^{-1} \bP^T \bV = \bV$, so that $(\bI - \cP)\bV = {\bf 0}$.
Therefore
\[
  \| \bA - \cP\bA \|  = \|(\bI - \cP) \bA\| 
  = \|(\bI - \cP)(\bI - \bV \bV^T) \bA\| \le 
   \|(\bI - \cP)\| \|(\bI - \bV \bV^T) \bA\|.
\]
It is well known that 
\[
   \|\bI - \cP\| = \| \cP \|  = \|(\bP^T \bV)^{-1} \| 
\]
so long as $\cP \ne {\bf 0} \ {\rm or} \ \bI$; see, e.g.,~\cite{Szy06}.
This establishes the first result.  The second follows from the fact that
\[
    \|(\bI - \bV \bV^T) \bA\| = \| \bA - \bV \bS \bW^T\| = \sigma_{k+1} 
\]
when $\bV$ consists of the leading $k$ left singular vectors of $\bA$.
\sq

\medskip
Now let $\bV \bS \bW^T \approx \bA$ be a rank-$k$ SVD of $\bA$.   
(The singular vectors play a crucial role in this analysis, even if $\bp$ and $\bq$ were selected using some scheme that did not reference them.)
In addition to the interpolatory projector 
$\cP = \bV(\bP^T\bV)^{-1}\bP^T$ that operates 
on the left of $\bA$, we shall also use 
$\cQ = \bQ(\bW^T \bQ)^{-1} \bW^T$, which operates
on the right of $\bA$.\ \ 
Assuming that $\bP^T\bV$ and $\bW^T\bQ$ are invertible, define the error constants
\[ \eta_p \equiv \|(\bP^T\bV)^{-1}\|, \qquad
   \eta_q \equiv \|(\bW^T\bQ)^{-1}\|. \]
Lemma~\ref{DEIM_bound} implies
\begin{equation} \label{eq:iebnd}
   \| \bA (\bI - \cQ) \| \le \eta_q @@\sigma_{k+1}
   \quad {\rm and } \quad
      \| (\bI - \cP) \bA \| \le \eta_p @@\sigma_{k+1}.
\end{equation}
The next lemma shows that these bounds on the
error of the interpolatory projection of $\bA$ onto the select columns and rows also apply to the \emph{orthogonal} projections of $\bA$ onto the same column and row spaces.

\begin{lemma}
\label{ProjErr}
Suppose the row and column indices $\bp$ and $\bq$ give
full rank matrices $\bC = \bA(@@:@@,\bq) = \bA\bQ\in\IR^{m\times k}$ and $\bR = \bA(\bp,@@:@@)=\bP\bA\in\IR^{k\times n}$, 
with finite error constants $\eta_p$ and $\eta_q$, and suppose that $k<\min\{m,n\}$.  Then
\[
   \| (\bI - \bC \bC^I) \bA \| \le \eta_q @@\sigma_{k+1}
   \quad {\rm and}\quad
   \| \bA(\bI - \bR^I \bR) \| \le \eta_p @@\sigma_{k+1}.
\]
\end{lemma}
{\it Proof:}
Using the formula $\bC= \bA\bQ$, we have 
$\bC^I = (\bC^T\bC)^{-1}\bC^T = (\bQ^T\bA^T\bA\bQ)^{-1}(\bA\bQ)^T$, 
so the orthogonal projection of $\bA$ onto $\Ran(\bC)$ is
\[
    \bC \bC^I \bA  = (\bA \bQ(\bQ^T \bA^T \bA \bQ)^{-1} \bQ^T \bA^T)\bA
      = \bA (\bQ (\bQ^T \bA^T \bA \bQ)^{-1} \bQ^T \bA^T\bA).
\]
Hence the error in the orthogonal projection of $\bA$ is
\begin{equation} \label{eq:lemid1}
    (\bI - \bC \bC^I) \bA = \bA(\bI - \bPhi), \ \ {\rm where} \ \ \bPhi =  \bQ(\bQ^T \bA^T \bA \bQ)^{-1} \bQ^T \bA^T\bA.
\end{equation}
Note that $\bPhi$ is an oblique projector onto $\Ran(\bQ)$, so  $\bPhi \bQ = \bQ$.  
Therefore, $\bPhi \cQ = \cQ$, since
\[
\bPhi \cQ =  \bPhi \bQ (\bW^T\bQ)^{-1} \bW^T =  \bQ (\bW^T\bQ)^{-1} \bW^T = \cQ.
\]
This implies that
\[
  \bA(\bI - \bPhi) = \bA(\bI - \bPhi)(\bI - \cQ) =  (\bI - \bC \bC^I) \bA (\bI - \cQ),
\]
and so from~(\ref{eq:lemid1}) we have
\begin{eqnarray*}
\|(\bI - \bC \bC^I) \bA \| &=& \| \bA(\bI - \bPhi)\| \\
     &=& \|  (\bI - \bC \bC^I) \bA (\bI - \cQ)\| \\
     &\le& \|  \bI - \bC \bC^I \| \|\bA (\bI - \cQ)\| \\
     &\le& \eta_q@@\sigma_{k+1}.
\end{eqnarray*}
The last line follows from the bound~(\ref{eq:iebnd}) and the fact that $ \|  \bI - \bC \bC^I \| = 1$, since $\bC\bC^I$ is an orthogonal projector and $k<\min\{m,n\}$.

A similar argument shows that 
\[
  \bA(\bI - \bR^I \bR) = (\bI - \bPsi) \bA  
\]
where $\bPsi = \bA \bA^T \bP (\bP^T\bA\bA^T\bP)^{-1}\bP^T$, and also that
\[
  (\bI - \bPsi)\bA = (\bI - \cP)(\bI - \bPsi)\bA = 
(\bI - \cP) \bA (\bI - \bR^I \bR),
\]
from which follows the error bound
\[
  \|\bA(\bI - \bR^I \bR)\|   \le \|(\bI - \cP) \bA \| \|\bI - \bR^I \bR\| \le \eta_p @@\sigma_{k+1}.\qquad \sq
\]

\medskip
The main result on approximation of $\bA$ by $\bC\bU\bR$
readily follows from combining this last lemma with a basic CUR analysis technique used by Mahoney and Drineas~\cite[eq.~(6)]{MD09}.
 
\begin{theorem}
\label{basic_bound}
Given $\bA\in\IR^{m\times n}$ and $1\le k<\min\{m,n\}$, let $\bC = \bA(@@:@@,\bq)\in\IR^{m\times k}$ and $\bR = \bA(\bp,@@:@@)\in\IR^{k\times n}$ 
with finite error constants $\eta_p$ and $\eta_q$, and set 
$\bU = \bC^I\bA\bR^I$.  Then
\[
   \| \bA - \bC \bU \bR \| \le (\eta_p + \eta_q)@@ \sigma_{k+1}.
\]
\end{theorem}
\noindent{\it Proof:} 
From the definitions,
\[
  \bA - \bC \bU \bR = \bA - \bC \bC^I \bA \bR^I \bR 
                    = (\bI - \bC \bC^I )\bA + \bC\bC^I \bA (\bI - \bR^I \bR).
\]
Applying Lemma~\ref{ProjErr},
\begin{eqnarray*}
  \|\bA - \bC \bU \bR \| 
                    &\le& \|(\bI - \bC \bC^I )\bA \| +  \|\bC\bC^I \| \| \bA (\bI - \bR^I \bR) \| \\
                    &\le& \eta_q@@\sigma_{k+1} +  \eta_p@@ \sigma_{k+1}\\
                    &=& ( \eta_p  +  \eta_q )@@\sigma_{k+1},
\end{eqnarray*}
since $\|\bC\bC^I\| = 1$.  \sq

\medskip
Theorem~\ref{basic_bound} shows that $\bC\bU\bR$ is within 
a factor of $\eta_p +\eta_q$ of the optimal rank-$k$ 
approximation, hence these error constants suggest a 
way to assess a wide variety of column/row selection 
schemes. 
The quality of the approximation is controlled by the conditioning of the selected $k$ rows of the dominant $k$ (exact) singular vectors.
If those singular vectors are available as part of the
column/row selection process, then Theorem~\ref{basic_bound}
provides an \emph{a posteriori} bound requiring only the fast ($\cO(k^3)$) computation of $\eta_p$ and $\eta_q$, and thus
could suggest methods for adjusting either $k$ or the point selection process to reduce the error constants.  In this context, notice that if $\bV\bS\bW^T$ is only an \emph{approximation} to the optimal rank-$k$ SVD with $\bV$ and $\bW$ having orthonormal columns (as computed, for example, using the incremental QR algorithm described in the next section), the preceding analysis gives
\begin{eqnarray}
\|\bA - \bC\bU\bR\| 
   \!\!&\le&\!\! \|(\bI-\bC\bC^I)\bA\| + \|\bA(\bI-\bR^I\bR)\| \nonumber \\
   \!\!&=&\!\! \|\bA(\bI-\cQ)\| + \|(\bI-\cP)\bA\| \nonumber \\
   \!\!&\le&\!\! \|(\bW^T\bQ)^{-1}\| \|\bA(\bI-\bW\bW^T)\| + \|(\bP^T\bV)^{-1}\| \|(\bI-\bV\bV^T)\bA\|, \label{eq:approxbnd}
\end{eqnarray}
showing how $\sigma_{k+1}$ in Theorem~\ref{basic_bound} is replaced by the error in the approximate SVD through $\|\bA(\bI-\bW\bW^T)\|$ and $\|(\bI-\bV\bV^T)\bA\|$. 
In this case $\|(\bW^T\bQ)^{-1}\|$ and $\|(\bP^T\bV)^{-1}\|$ are computed using the approximate singular vectors in $\bV$ and $\bW$, rather than the exact singular vectors in the theorem.
Alternatively, if one has probabilistic bounds for $\eta_p$ and $\eta_q$, 
then Theorem~\ref{basic_bound}  immediately gives a probabilistic bound for $\| \bA - \bC \bU \bR \|$.

Numerical examples in Section~\ref{sec:examples} compare how the error constants evolve as $k$ increases for the DEIM-CUR 
factorization and several other factorizations based on leverage scores.  

\subsection{Interpretation of the bound for DEIM-CUR}

For DEIM-CUR, we can ensure  the hypotheses of Theorem~\ref{basic_bound} are satisfied and bound the error constants.  Suppose the DEIM points are selected using the exact rank-$k$ SVD $\bA \approx \bV\bS\bW^T$.\ \ Lemma~\ref{DEIM_inv} ensures that the matrices
$\bP^T\bV$ and $\bW^T\bQ$ are invertible, so $\eta_p$ and $\eta_q$ are finite.  The DEIM strategy also gives full rank $\bC$ and $\bR$ matrices, presuming $k\le {\rm rank}(\bA)$.  To see this, note that for any unit vector $\by\in\IR^k$,
\[
   \bC\by =  \bA \bQ@\by = \bV \bS \bW^T \bQ@\by + \bE\bQ@ \by,
\]
where $\bE = \bA - \bV\bS\bW^T$.\ \ 
Since $\bV^T \bE = {\bf 0},$
\[
   \| \bC\by\|^2  = \|  \bA \bQ@\by\|^2  =\|  \bV \bS \bW^T \bQ\by\|^2  + \| \bE\bQ@\by \|^2.
\]
Since  $\| \bW^T \bQ\by\| \ge \|\by\|/\|(\bW^T\bQ)^{-1}\| = 1/\eta_q$,
\[
   \| \bC\by\|  \ge \|  \bV \bS \bW^T \bQ\by\| \ge \sigma_k / \eta_q > 0.
\]
Thus $\bC$ must be full rank.   A similar argument shows $\bR$ to be full rank as well.
\medskip

The examples in Section~\ref{sec:examples} illustrate that $\eta_p$ and $\eta_q$ are often quite modest for the DEIM-CUR approach, e.g., $\cO(100)$.\ \  However, worst-case  bounds permit significant growth in $k$ that is generally not observed in practice.  
We begin by stating a bound on this growth developed by
Chaturantabut and Sorensen~\cite[Lemma~3.2]{CS10}.

\begin{lemma} 
For the DEIM selection scheme derived above,
\[\eta_p \le {(1+\sqrt{2m})^{k-1}\over \|\bv_1\|_\infty}, \qquad
   \eta_q \le {(1+\sqrt{2n})^{k-1} \over \|\bw_1\|_\infty},\]
where $\bv_1$ and $\bw_1$ denote the first columns of $\bV$ and $\bW$.
\end{lemma}

Motivated by recent work by Drma\v{c} and Gugercin~\cite{DG} on a modified DEIM-like algorithm for model reduction, we can improve this bound considerably.
\begin{lemma} \label{lemma:new_eta_bound}
For the DEIM selection scheme derived above,
\[
   \eta_p < \sqrt{mk\over 3} \ {2^k}, \qquad
   \eta_q < \sqrt{nk\over 3} \ {2^k}.
\]
\end{lemma}
\noindent{\it Proof:} 
We shall prove the result for $\eta_p$; the result for $\eta_q$ follows similarly.
As usual, let $\bV \in \IR^{m \times k}$ have orthonormal columns, 
and let $\bp = {\rm DEIM}(\bV)$ denote the row index vector derived 
from the DEIM selection scheme described above.  
Let $\bP = \bI(\,:\,,\bp)$ so that $\bP^T \bV = \bV(\bp,\,:\,)$.  

Without loss of generality, assume the DEIM index selection gives 
$\bp = [1, 2, \ldots, k]^T$.  (Otherwise, introduce a permutation matrix to the argument that follows.)
As described in section~\ref{sec:deim}, the DEIM index selection 
is precisely the index selection of LU decomposition with partial pivoting,
so one can write
\[
    \bV = \bL\kern1pt \bT,
\]
where the nonsingular matrix $\bT\in\IR^{k\times k}$ is upper triangular 
and $\bL\in\IR^{m\times k}$ is unit lower triangular with $ | \bL(i,j) | \le 1$, 
$\bL(j,j) = 1, \ 1 \le j \le k$ and $\bL(i,j) = 0, \ j > i $.

Let $\bL_1 \equiv \bL(1:k,1:k)$.  Then $\bV(\bp,\,:\,) = \bL_1 \bT$ and thus
\[
 \eta_p \equiv \| (\bP^T \bV)^{-1} \| = \|(\bL_1 \bT)^{-1} \| \le \| \bT^{-1}\| \| \bL_1^{-1}\|.
\]
(The linear independence of the columns of $\bV$ ensure that $\bL_1$ and $\bT$ are invertible.)
Upper bounds for $\|\bT^{-1}\|$ and $\| \bL_1^{-1}\|$ will give an upper bound for $\eta_p$.

To bound $\| \bT^{-1}\|$, let $\by\in\IR^k$ be a unit vector such that 
$\| \bT^{-1} \by \| = \| \bT^{-1}\|$.
Then
\[
   \| \bT^{-1} \| = \| \bT^{-1} \by \| = \| \bV \bT^{-1} \by \| = \| \bL \by \|.
\]
Now
\[
   \| \bL \by \| \le \sqrt{m}\, \|\bL\by\|_\infty = \sqrt{m}\, |\bfe_j^T \bL \by |,
\]
for some index $j\in\{1,\ldots, m\}$.\ \ 
By the Cauchy--Schwarz inequality and the bound $|\bL(i,j)|\le 1$, 
\[
  |\bfe_j^T \bL \by | \le \|\bL^T \bfe_j\| \| \by \| \le \sqrt{k}\cdot 1,
\]
and so it follows that $\| \bT^{-1} \| \le \sqrt{mk}$.

The inverse of $\bL_1$ can be bounded using forward substitution.
Let $\bL_1 \bz = \by$, where $\|\by\| = 1$ and $\| \bz \| =  \| \bL_1^{-1}\|$.
Forward substitution provides
\begin{eqnarray*}
   \zeta_1 &=& \gamma_1 \\
   \zeta_i &=& \gamma_i - \sum_{j=1}^{i-1} \lambda_{ij} \zeta_j, \rlap{\qquad $i = 2,\ldots, k$,}\\
\end{eqnarray*}
where $\zeta_i = \bz(i)$, $\gamma_i = \by(i)$ and $\lambda_{ij} = \bL(i,j)$.
We now use induction to prove
\[
    | \zeta_i | \le 2^{i-1}, \rlap{\qquad $1 \le i \le k$.}
\]
First note that $\zeta_1 = \gamma_1$, so $|\zeta_1| \le |\gamma_1| \le 1 = 2^0$
to establish the base case.  Assume for some $i \ge 1$ that
\[
    | \zeta_j | \le 2^{j-1}, \rlap{\qquad $1 \le j \le i$.}
\]
Then
\begin{eqnarray*}
   | \zeta_{i+1} | = \bigg|  \gamma_{i+1} - \sum_{j=1}^{i} \lambda_{ij} \zeta_j \bigg|
                   &\le& |\gamma_{i+1}| + \sum_{j=1}^{i} |\lambda_{ij}| |\zeta_j | \\
                   &\le&  1 + \sum_{j=1}^{i} 1 \cdot 2^{j-1}
                   =  1 + \sum_{j=0}^{i-1}  2^j 
                   = 1 + (2^i - 1) = 2^i
\end{eqnarray*}
to complete the induction.
Now since $\|\bz\| = \|\bL_1^{-1}\|$,
\[
   \| \bL_1^{-1} \|^2 = \bz^T \bz = \sum_{i=1}^k |\zeta_i|^2 \le \sum_{i=0}^{k-1} 4^i = (4^k - 1)/3.
\]
Thus $ \| \bL_1^{-1} \| < 2^k/\sqrt{3}$, which, together with the bound on the inverse of $\bT$, provides
the final result for $m>k$:
\[
     \eta_p \equiv \| (\bP^T \bV)^{-1} \| < \sqrt{mk\over 3}\, 2^k.
\]
If $m=k$, then $\eta_p=1$, and the result holds trivially.
\quad \sq

\medskip
Note that this proof only relies on the orthonormality of the columns of $\bV$ and $\bW$,
and hence it applies when the DEIM selection scheme is applied to approximate 
singular vectors, as in CUR error bound in~(\ref{eq:approxbnd}).

Lemma~\ref{lemma:new_eta_bound}  
was inspired by the proof technique developed by 
Drma\v{c} and Gugercin~\cite{DG} to bound $\|(\bP^T \bV)^{-1} \|$,
when $\bP$ is selected by applying a pivoted rank-revealing 
QR factorization scheme to $\bV$.\ \   
Note that this new bound is on the same order of magnitude as the 
Drma\v{c}--Gugercin scheme.  In practice, their scheme
seems to give slightly smaller growth that is more consistent over a wide range of examples.
Neither scheme experienced exponential growth over very extensive testing.
For the DEIM approach, this absence of exponential growth is closely related
to decades of experience with Gaussian elimination with partial pivoting.   
Element growth in $\bT$ is bounded by a factor of $2^{k-1}$ 
(for a $k \times k$ matrix), and there is an example that achieves this 
growth.   Nevertheless, this algorithm is almost exclusively used to solve
linear systems because such growth is never experienced.%
\footnote{See, for example, the extensive numerical tests 
involving random matrices described 
in~\cite[lecture~22]{TB97} and~\cite{TS90}.
Interestingly, in the experiments of Trefethen and Schreiber~\cite{TS90},
random matrices with orthonormal columns tend to have slightly larger 
growth factors than Gaussian matrices, though both cases are very 
far indeed from the exponential upper bound.}
Indeed, a similar near worst case example can be constructed for DEIM,
although this growth has not been observed in practice.

\medskip
\noindent{\bf A Growth Example:}
We now construct an orthonormal matrix $\bV$ with the property
\begin{equation}\label{eq:growthex}
   {1\over \sqrt{8}}\ 2^{k}  <  \eta_p \equiv \|(\bP^T \bV)^{-1} \| < \sqrt{m k\over 3}\ 2^k
\end{equation}
where $\bP^T \bV = \bV(\bp,\,:\,)$ with $\bp = {\rm DEIM}(\bV)$.  
To construct $\bV$, begin by defining
\[ \bL := \left[\begin{array}{cccc}
1 & & & \\
-1 & 1 & & \\
\vdots & \ddots & \ddots & \\
-1 & \cdots & -1 & 1 \\
-1 & \cdots & -1 & -1 \\
\vdots &     &   &  \vdots \\
-1 & \cdots & -1 & -1
           \end{array}\right]\in \IR^{m\times k}.\]
Now construct $\bV\bT_1 \equiv \bL$ as an economy-sized QR factorization of
$\bL$ (with no column pivoting).  Since the columns of $\bL$ are
linearly independent by construction, $\bT_1\in\IR^{k\times k}$ is invertible;
define $\bT \equiv \bT_1^{-1}$, so that $\bV = \bL\kern1pt\bT$.
(Note that $\bT$ plays the same role it does in the proof of Lemma~\ref{lemma:new_eta_bound}.)
If the DEIM procedure is applied to $\bV$, then by construction
$\bp = [1, 2, \ldots, k]$  (in exact arithmetic):
during the DEIM procedure,  the relations
\[
     \bell_j\kern1pt\tau_{jj} = \bv_j -  \bV_{j-1}^{} ( \bP_{j-1}^T \bV_{j-1}^{})^{-1} \bP_{j-1}^T \bv_j, \rlap{\qquad $j > 1$}
\]
hold, with $\bell_j  = \bL(\,:\,,j)$, $\tau_{jj} = \bT(j,j)$, 
$\bv_j = \bV(\,:\,,j)$ , $\bP_{j-1} = \bI(\,:\,,1:j-1)$
and $\bV_{j-1} = \bV(\,:\,,1:j-1)$.  Thus $\bp(j) = j , \ j>1 $ and it is easily seen that $\bp(1) = 1$.

Note that $ \bV \bT^{-1} = \bL $ implies $\bT^{-T} \bT^{-1} = \bL^T \bL$ hence $\| \bT \| = 1/\sigma_k$,
where $\sigma_k $ is the smallest singular value of $\bL$.\ \ 
Let $\by$ be the corresponding right singular vector, so that
\[
   \sigma_k^2 = \by^T \bL^T \bL \by.
\]
We claim that $\sigma_k \ge \sqrt{2}$.\ \ 
To see this, write $\bL$ in the form
\[
   \bL = \left[ \begin{array}{c}
                     \bI_k \\
                     {\bf 0}
                \end{array}
           \right]
         - \left[ \begin{array}{c}
                     \bL_0 \\
                      \bE
                \end{array}
           \right],
\]
where 
\[ \bI_k = \left[\begin{array}{cccc}
                   1 \\ & 1 \\ & & \ddots \\ & & & 1
            \end{array}\right], 
    \quad
    \bL_0 = \left[\begin{array}{cccc}
               0 \\ 1 & \ddots \\ \vdots & \ddots & \ddots \\
               1 & \cdots & 1 & 0\end{array}\right],
     \quad 
     \bE = \left[\begin{array}{cccc}
               1 & 1 & \cdots & 1 \\ 
               1 & 1 & \cdots & 1 \\ 
               \vdots & \vdots & \ddots & \vdots \\
               1 & 1 & \cdots  & 1\end{array}\right]
            = \bff\kern1pt \bfe^T,
\]
for $\bff = [1,\ldots, 1]^T \in \IR^{m-k}$.\ \ 
\begin{eqnarray*} 
 \bL^T \bL &=& \bI_k - \bL_0 - \bL_0^T + \bL_0^T \bL_0 + \bE^T \bE \\[.25em]
           &=& \bI_k - (\bfe \bfe^T - \bI_k) + \bL_0^T\bL_0 + (m-k)\bfe\bfe^T \\[.25em]
           &=& 2\kern1pt \bI_k + \bL_0^T \bL_0 + (m - k -1) \bfe \bfe^T \\[.25em]
           &=& 2\kern1pt \bI_k +  \bM,
\end{eqnarray*}
where $\bM:= \bL_0^T \bL_0 + (m - k -1) \bfe \bfe^T$ 
is symmetric and positive semidefinite whenever $m>k$.
Thus
\[
   \sigma_k^2 = \by^T \bL^T \bL \by = \by^T (2\bI_k +  \bM) \by \ge 2
\]
and hence it follows that
\[
   \| \bT \| \le 1/\sqrt{2}.
\]
This implies
\[
   \|(\bP^T \bV)^{-1} \| = \| \bT^{-1} \bL_1^{-1} \| \ge \frac{\| \bL_1^{-1} \|}{\| \bT \|} \ge
     \sqrt{2}\ \|\bL_1^{-1} \|.
\]
To complete the lower bound, we must analyze $\|\bL_1^{-1}\|$.  Forward substitution gives $\bL_1^{-1} \bfe_1 = [1, 1, 2, 4, \ldots , 2^{k-2}]^T$ and thus
\[
    \|\bL_1^{-1} \| > \|\bL_1^{-1} \bfe_1\| = \sqrt{1 + (4^{k-1} - 1)/3} \ > 2^{k-2}.
\]
We arrive at the lower bound
\[
   \eta_p \equiv \|(\bP^T \bV)^{-1} \| 
     \ge  \sqrt{2}\ \|\bL_1^{-1} \| \  > \sqrt{2}\cdot 2^{k-2},
\]
and thus for this choice of $\bV$, the DEIM error constant satisfies
\[
  {1\over\sqrt{8}}\ 2^k  < \eta_p < \sqrt{{m k\over 3}} \ 2^k.
\]

This example is interesting because it relies on the behavior of the
classic example for growth in LU decomposition~\cite[lecture~22]{TB97}.
However, in this case the pathological growth is caused by $\bL$ 
and not by $\bT$.

\section{Incremental QR Factorization} \label{sec:qr}

The DEIM point selection process presumes access to the first $k$ left and right singular vectors of $\bA\in\IR^{m\times n}$.  If either $m$ or $n$ is of modest size (say ${}\le 1000$) and $\bA$ can be stored as a dense matrix, \ library software for computing the ``economy sized'' SVD, e.g., \verb|[V,S,W] = svd(A,'econ')| in MATLAB, usually performs very well. 
For larger scale problems, the leading $k$ singular vectors can be computed using iterative methods, such as the Krylov subspace-based ARPACK software~\cite{LSY98} (used by MATLAB's {\tt svds} command), PROPACK~\cite{Lar05}, IRLBA~\cite{BR13}, or the Jacobi--Davidson algorithm~\cite{Hoc01}.  Randomized SVD algorithms provide an appealing alternative with probabilistic error bounds~\cite{HMT11}. 
Here we describe another approach that satisfies a deterministic error bound (Lemma~\ref{lem:incqr}) and \emph{only requires one pass through the matrix $\bA$}, a key property for massive data sets that cannot easily be stored in memory.

\begin{algorithm}[h!]
\begin{center}
\fbox{\begin{tabular}{p{10cm}}
\ \\[-25pt]
\begin{tabbing}
defghasijkl\=bbbb\=cccc\= \kill
$ \ \ $ {\bf Input:} \> $ \bA$,  an $m \times n$ matrix \\
$ \ \ $ \> {\it tol}, a positive scalar controlling the accuracy of the factorization\\[5pt]
$\  $  {\bf Output:}  \> $\QQ$, an $m \times k$ matrix with orthonormal columns\\
$\ \ \ \ \ \ \ $ \> $\RR$, a $k \times n$ rectangular matrix  \\
$\ \ \ \ \ \ \ $ \>\>\>  with $\bA \approx  \QQ \RR$ \\[5pt]
$\ \ \ $ Choose $k \ll \min(m,n)$\\
$\ \ \ $ Compute the QR factorization $\bA(:,1:k) = \QQ \RR$, with $\QQ \in \IR^{n\times k}$ and $\RR\in\IR^{k\times m}$\\  
%
%
$ \ \ \ $  rownorms($i$) = $\|\RR(i,:)\|^2$ for $i = 1, \ldots, k$ \\
$\ \ \ $   $j = k+1$\\
$\ \ \ $ {\bf while} $j \le n$\\
$ \ \ \ $  \> $\bfa = \bA(@:@,j) ; \ \rr = \QQ^T\bfa; \ \  \bff = \bfa - \QQ \rr; \ \   \rho = \|\bff \|; \ \ \qq = \bff/\rho$\\
$\ \ \ $  \>   $\QQ = [\QQ, \ \qq];  \ \   \RR = \left[ \begin{array}{cc} \RR  &  \rr\\
                                                                {\bf 0} & \rho \\
                                               \end{array} \right]$ \\[5pt]
$\ \ \ $  \> rownorms($i$) = rownorms($i$) + $\rr(i)^2$ for $i = 1,\ldots, k$ \\                                               
$\ \ \ $  \> rownorms($k+1$) ${}= \rho^2; \ \  $ \\
$ \ \ \ $ \>  FnormR = {\rm sum}(rownorms); \\
$ \ \ \ $ \> $[\sigma, i_{\rm min}] = {\rm min}({\rm rownorms}(1:k+1));$ \\[5pt]
$ \ \ \ $ \> {\bf if} $\sigma > ({\it tol}^2)*({\rm FnormR} - {\rm rownorms}(i_{\rm min})$ \\
%
%
$ \ \ \ $ \>\> \% \emph{no deflation} \\
$ \ \ \ $ \>\> $ k = k + 1;$ \\
$ \ \ \ $ \> {\bf else} \% \emph{deflation required}\\
%
%
$ \ \ \ $ \>\> {\bf if} $i_{\rm min}  < k+1$\\
%
%
$ \ \ \ $ \>\>\> $\RR(i_{\rm min},@:@) = \RR(k+1,@:@); \ \QQ(@:@,i_{\rm min}) = \QQ(@:@, k+1)$ \\
$ \ \ \ $ \>\>\> $\ {\rm rownorms}(i_{min}) = {\rm rownorms}(k+1) $ \\
%
%
$ \ \ \ $ \>\> {\bf end} \\
%
%
$ \ \ \ $ \>\> \%  \emph{delete the minimum norm row of $\RR$}\\ 
$ \ \ \ $ \>\> $\QQ = \QQ(@:@,1:k); \RR = \RR(1:k,@:@) $ \\
$ \ \ \ $ \> {\bf end} \\
$ \ \ \ $ \> $j = j+1$\\
$ \ \ \ $  {\bf end} \\
\end{tabbing}\\[-23pt]\end{tabular}}
\end{center}

\vspace*{-10pt}
\caption{Incremental QR low rank approximate factorization}
\label{fig:QRincremental}
\end{algorithm}

\begin{figure}[t!]
\begin{center}
\begin{picture}(300,50)
\put(-40,-20){$\bA(@:@,1:j) = {}$}
\put(-25,-80){
\begin{minipage}[t]{2.5in}\begin{center} \hspace*{30pt}\includegraphics[scale=0.35]{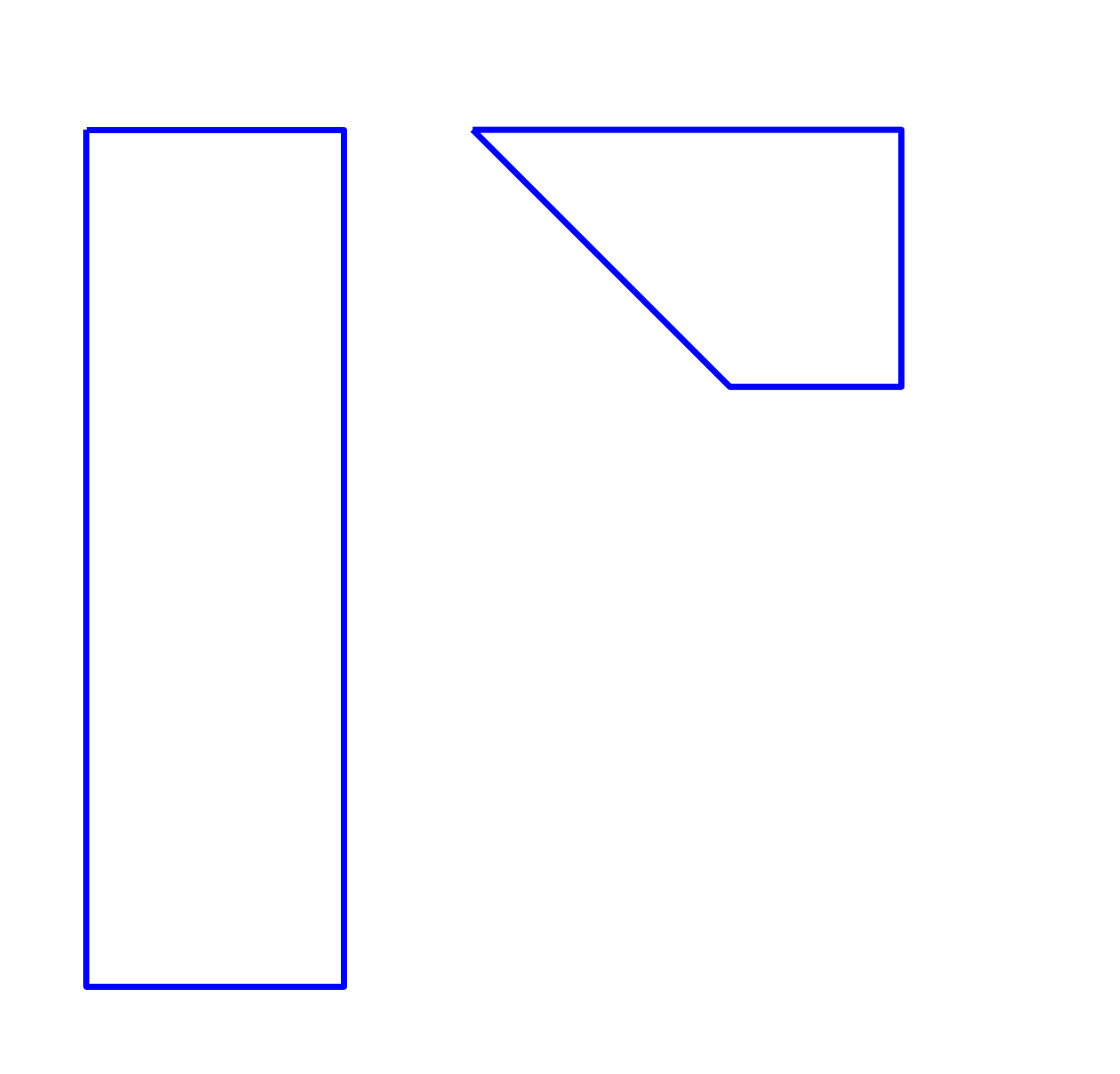}\\  Partial $\bQ\RR$ factorization\end{center} \end{minipage}}

\put(160,-20){\small $\bA(:,1\!:\!j\!+\!1) = {}$}
\put(177,-80){\begin{minipage}[t]{2.5in}\begin{center}\hspace*{30pt} \includegraphics[scale=0.35]{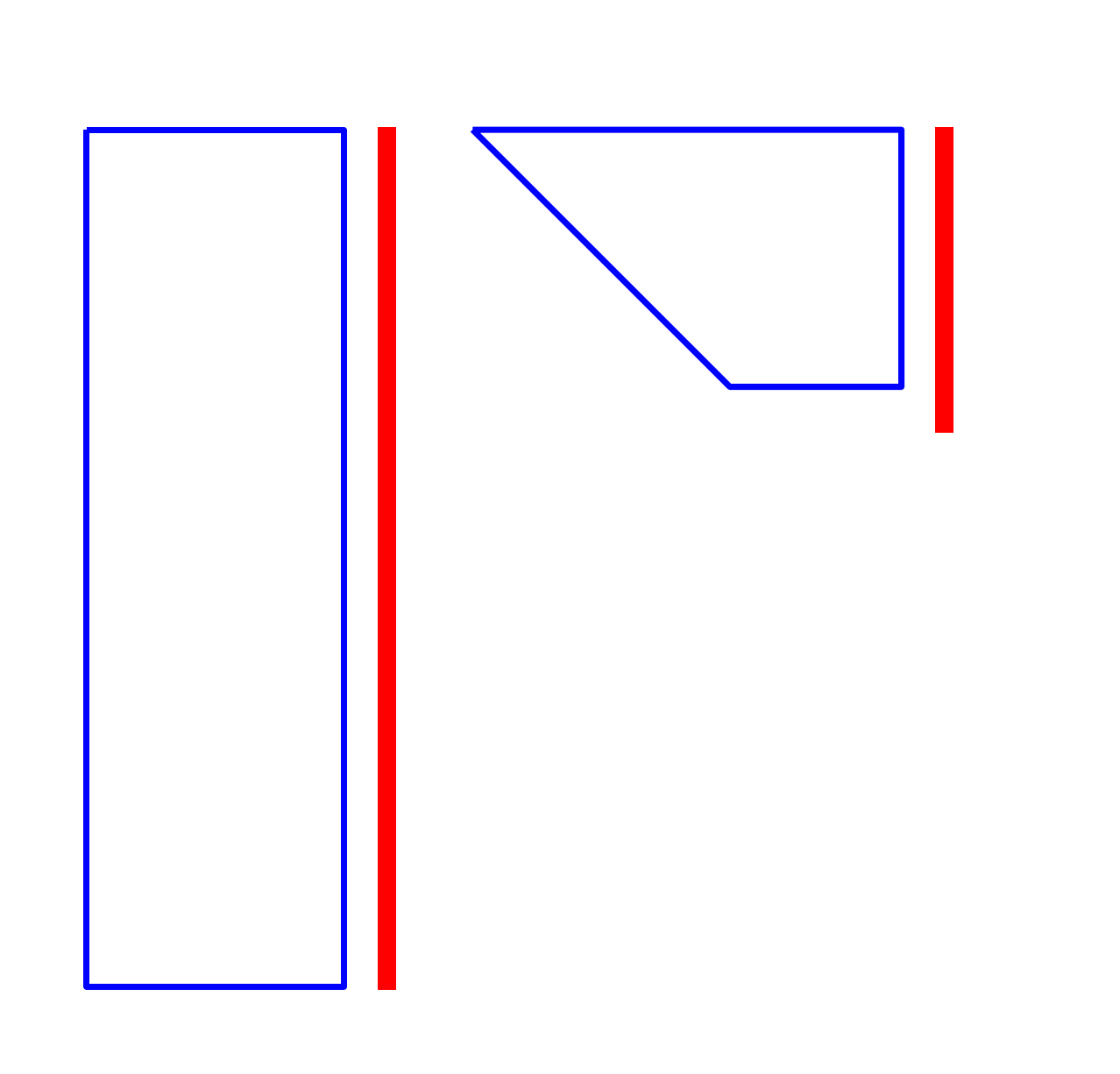}\\  Extend with Gram--Schmidt\end{center} \end{minipage}}

\put(-95,-230){\begin{minipage}[t]{3in}\begin{center} \includegraphics[scale=0.35]{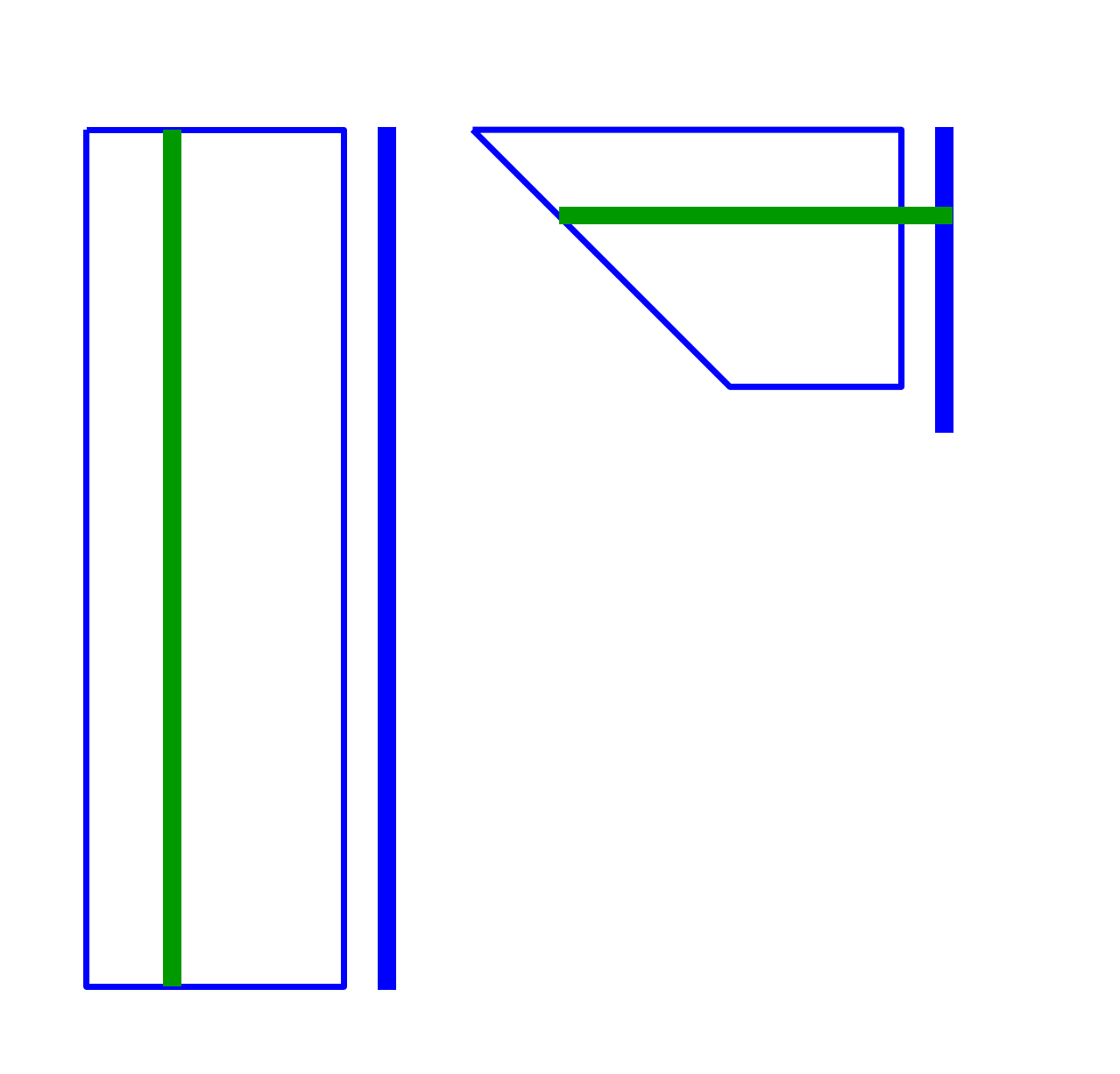}\\  Find $\qq_i$ with $\|\RR(i,:)\|^2 < \tol^{@2}\,(\|\RR\|_F^2 - \|\RR(i,:)\|^2)$\end{center} \end{minipage}}

\put(90,-230){\begin{minipage}[t]{2in}\begin{center} \includegraphics[scale=0.35]{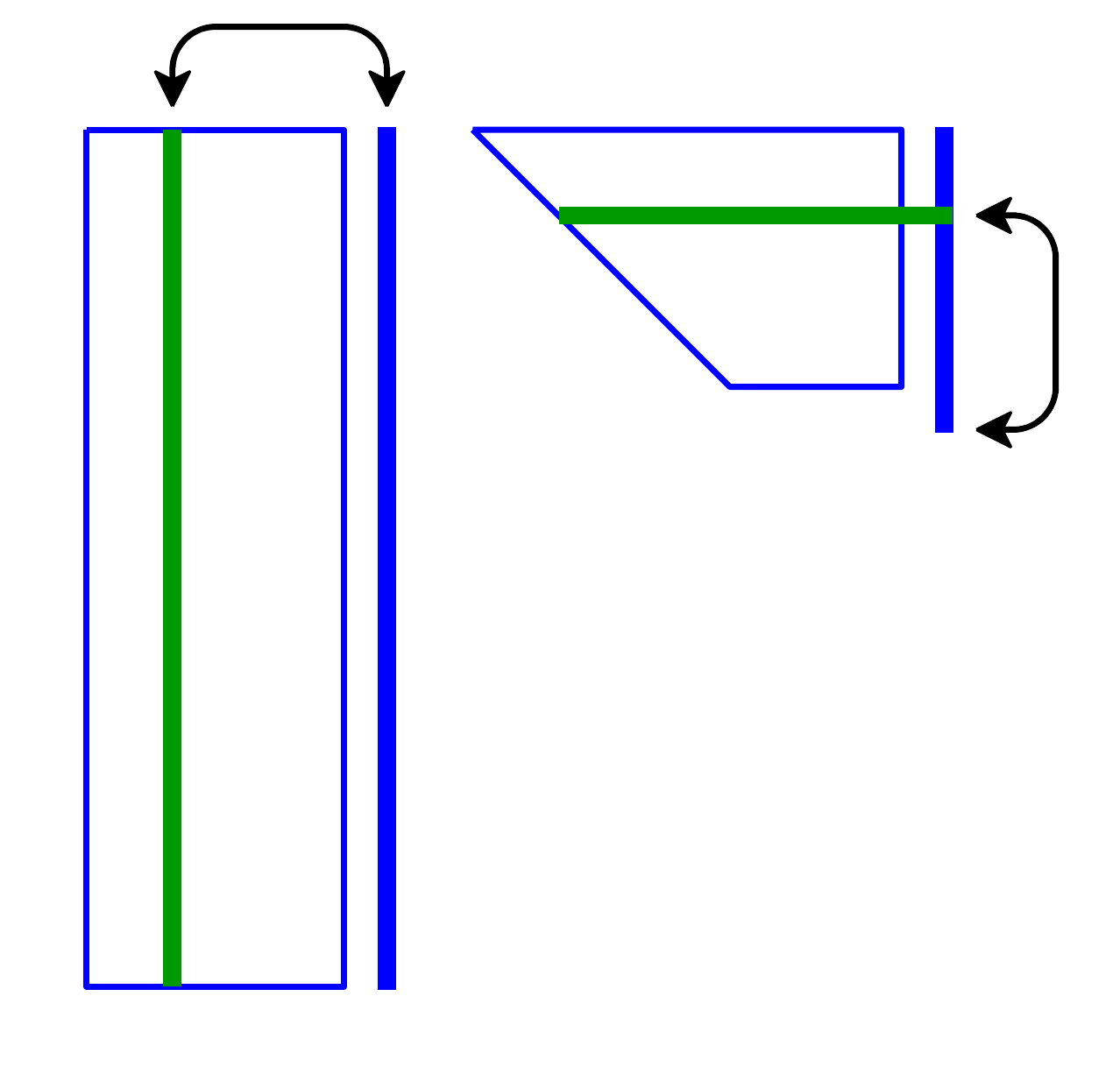}\\ Replace $\qq_i$, $\RR(i,:)$ \end{center} \end{minipage}}

\put(255,-230){\begin{minipage}[t]{2in}\begin{center} \includegraphics[scale=0.35]{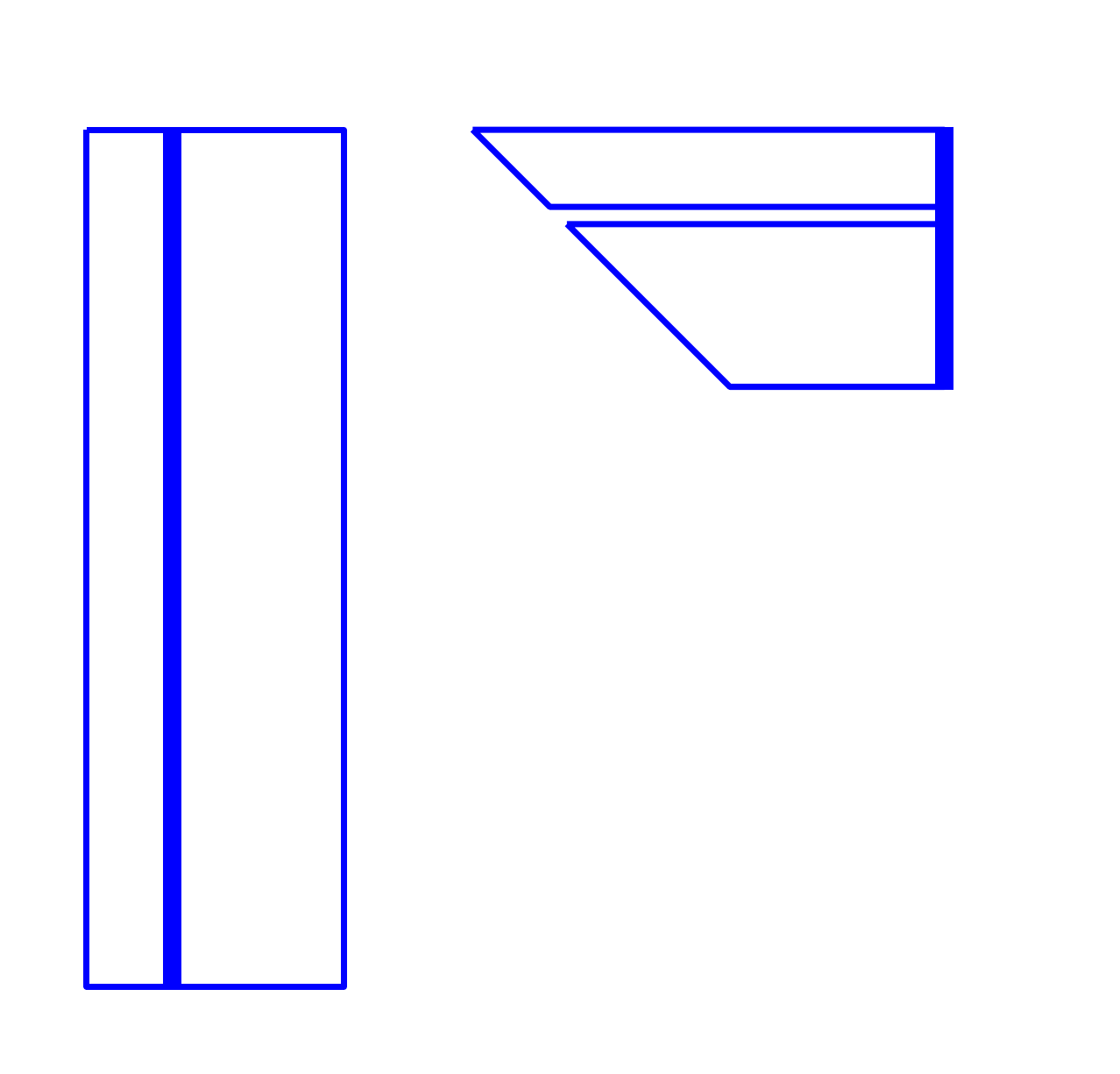}\\ Truncate last column\\ of $\QQ$ and row of $\RR$\end{center} \end{minipage}}
\end{picture}
\vspace*{3.35in}
\end{center}
\caption{\label{fig:QRupdate}
Diagram illustrating the QR update procedure.}
\end{figure}

This approach is based on an {\it incremental} low
rank $\bA\approx\bQ \RR $ approximation, where $\QQ\in\IR^{n\times k}$ has orthonormal columns and $\RR\in\IR^{k\times m}$ is upper triangular.  
(In this section only, $\QQ$ and $\RR$ denote different quantities from elsewhere in the paper.)
Take the dense (economy sized) SVD $\RR = \widehat{\bV} \bS \bW^T$, and put $\bV = \QQ \widehat{\bV} $ to get
\begin{equation} \label{eq:AQRSVD}
   \bA \approx  \QQ \RR = \bV \bS \bW^T.
\end{equation}
Incremental algorithms for building the QR factorization and SVD
have been proposed by Stewart~\cite{Ste99}, Baker, Gallivan, and Van Dooren~\cite{BGV12} and many others, as
surveyed in~\cite{Bak04}; these ideas are also closely related to rank-revealing QR factorizations~\cite{GE96}.  Algorithm~\ref{fig:QRupdate} differs from those of Stewart in its use of internal pivoting and threshold truncation in place of Stewart's column pivoting.   This distinction enables a
one-pass algorithm that is closely related to~\cite[Algorithm~1]{BGV12}.

The proposed method is presented in Algorithm~\ref{fig:QRincremental}, which proceeds at each step by orthogonalizing a column of $\bA$ against the previously orthogonalized columns. 
The rank of the resulting factors is controlled through an update-and-delete procedure that is illustrated in Figure~\ref{fig:QRupdate}.  
After orthogonalizing a column of $\bA$, 
the algorithm checks if any row of $\RR$ has small relative norm; 
if such a row exists, the corresponding column of $\bQ$ 
makes little contribution to the factorization, so that column of $\bQ$ and row of $\RR$ can be deleted at only a small loss of accuracy in the factorization.  
(Future columns of $\bA$ will not be orthogonalized
against the vector deleted from $\bQ$, so this direction can re-emerge if a later column in $\bA$ warrants it.)

Robust implementations of Algorithm~\ref{fig:QRupdate} should replace the classical Gram--Schmidt operations 
\[
   \rr = \QQ^T \bfa, \ \  \bff = \bfa - \QQ \rr
\]
with a re-orthogonalization step, as 
suggested by Daniel, Gragg, Kaufman, and Stewart~\cite{DGKS}:
\begin{eqnarray}
  \rr &=& \QQ^T \bfa \nonumber \\
  \bff &=& \bfa - \QQ \rr \nonumber \\
       &\ & \bc = \QQ^T \bff \label{eq:reorth1} \\
       &\ &  \bff = \bff - \QQ\bc \label{eq:reorth2}\\
       &\ &  \rr = \rr + \bc \label{eq:reorth3}\\
  \rho &=& \| \bff \| \nonumber \\
  \qq  &=&  \bff/\rho. \nonumber 
\end{eqnarray}
The extra steps~(\ref{eq:reorth1})--(\ref{eq:reorth3}) generally provide a $\QQ$ that is numerically orthogonal to working precision.
Pathological cases are easily overcome with some additional slight modifications;
see~\cite{GIR05} for a complete analysis.
Because this algorithm uses the classical Gram--Schmidt method, one can easily block it for parallel efficiency.

\subsection{Incremental QR Error Bounds}

At step $j$ the truncation criterion in Algorithm~\ref{fig:QRincremental} will delete row $\rr_i^T =  \bfe_i^T \RR_j \ $ if
\[
 \| \rr_i \| \le \tol\, \| \widehat{\RR}_j \|_F,
\]
where $\rr_i^T$ is the row of minimum norm and $\widehat{\RR}_j$ denotes $\RR_j$ with the $i$th row deleted.  This strategy has a straightforward error analysis, 
which, in light of the approximation~(\ref{eq:AQRSVD}), also 
implies an error bound on the resulting SVD.

\begin{lemma}\label{lem:incqr}
Let $\RR_j$ be the triangular factor at step $j$ of 
Algorithm~\ref{fig:QRincremental}, and $\QQ_j $ the corresponding orthonormal columns in the approximate QR factorization $\bA_j \approx \QQ_j \RR_j$, where $\bA_j$ consists of
the first $j$ columns of $\bA$. Then
\[
 \| \bA_j - \QQ_j \RR_j \|_F \le \tol \cdot d_j \cdot \| \RR_j \|_F, 
 \]
 where $d_j$ is the number of coloumn/row deletions that have been made up to and including step $j$.
 (Note that $\QQ_j \in \IR^{m\times (j-d_j)}$, $\RR\in\IR^{(j-d_j)\times j}$, and $d_n = \min\{m,n\} - k$, where $k = \rank(\QQ_n\RR_n)$.)
\end{lemma}

Before proving this lemma, we note that it gives a bound on the error in the resulting
approximate SVD of $\bA$.  Suppose $d_n$ deletions are made when this algorithm computes
the approximate factorization $\bA \approx \bQ\bR$ with tolerance $\tol$. 
Given the SVD $\bR = \widehat{\bV}\bS\bW^*$, set $\bV \equiv \bQ\widehat{\bV}$.
Then
\[ \|\bA - \bV\bS\bW^*\|_F \le \tol \cdot d_n \cdot  \|\RR\|_F.\]
{\it Proof of Lemma~\ref{lem:incqr}:}
The proof shall be by induction.  Let $\bE_j = \bA_j - \QQ_j  \RR_j$ and assume
\begin{equation} \label{eq:inductive}
\| \bE_j \|_F \le \tol \cdot d_j \cdot \| \RR_j \|_F.
\end{equation}
Orthogonalize column $j+1$ of $\bA$ using Gram--Schmidt to obtain
\[
 \bA_{j+1} = \QQ_{j+1}\RR_{j+1} + [\bE_j, {\bf 0}].
\]
If no deflation occurs at this step, the bound holds trivially
since 
\[
   \| \bE_{j+1} \|_F = \|[\bE_j, {\bf 0}]\|_F \le tol\cdot d_j \cdot \|\RR_j\|_F
   \le \tol \cdot d_{j+1} \cdot \| \RR_{j+1} \|_F,
\]
because $d_{j+1}= d_j$ and $\|\RR_j \|_F \le \| \RR_{j+1}\|_F$.

Suppose $\RR_j$ has dimension $k \times j$ (i.e., $k=j-d_j$).
Let $i$ be the index of the row of minimum norm and let $\widehat{\RR}_{j+1}$
be obtained by deleting the $i$th row of $\RR_{j+1}$.\ \ 
If $\rr_i^T  = \bfe_i^T \RR_{j+1} $ satisfies $\|\rr_i^T \| \le tol\cdot \| \widehat{\RR}_{j+1} \|_F$
then deflation occurs.
Deleting column $i$ of $\QQ_{j+1}$ and row $i$ of $\RR_{j+1}$
replaces $\QQ_{j+1}$ and $\RR_{j+1}$ with $\widehat{\QQ}_{j+1}$ and $\widehat{\RR}_{j+1}$.  Then
\[
\widehat{\QQ}_{j+1} \widehat{\RR}_{j+1} = \QQ_{j+1} ( \RR_{j+1} - \bfe_i^{} \rr_i^T),
\]
and 
\[
  \bA_{j+1} = \QQ_{j+1} (\RR_{j+1} - \bfe_i^{} \rr_i^T) + [\bE_j, {\bf 0}] + \QQ_{j+1}\bfe_i^{} \rr_i^T.
\]
Hence the deletion gives the overall error
\[
  \bE_{j+1} = \bA_{j+1} - \widehat{\QQ}_{j+1} \widehat{\RR}_{j+1} = [\bE_j, {\bf 0}] + 
   \QQ_{j+1}\bfe_i^{}\rr_i^T.
\]
Therefore, when $i < k+1$, the inductive assumption~(\ref{eq:inductive}) implies
\[
   \| \bE_{j+1}\|_F \le \| \bE_j \|_F + \|\rr_i^T \| \le \tol \cdot ( d_j \cdot \| \RR_j \|_F + \|\widehat{\RR}_{j+1} \|_F) \le \tol \cdot (d_j+1) \cdot \| \widehat{\RR}_{j+1}\|_F,
\]
since $\widehat{\RR}_{j+1}$ contains row $k+1$ of $\RR_{j+1}$, 
which must have a norm larger than the row marked for deletion.   
Since row $k+1$ of $\widehat{\RR}_{j+1}$ consists of just one
nonzero element, 
\[
 \| \widehat{\RR}_{j+1} \|^2_F \ge \| \widehat{\RR}_j \|_F^2 + \rho_{k+1,j+1}^2 \ge \| \RR_j \|_F^2,
 \]
 where $\rho_{k+1,j+1}$ is the element  $\RR_{j+1}(k+1,j+1) $  and $\widehat{\RR}_j$
 is the matrix $\RR_j$ with $i$th row deleted.
 If $i=k+1$, then the last row of $\RR_{j+1} $ is deleted and the desired inequality must hold, since $\RR_j$ is a submatrix of $\widehat{\RR}_{j+1}$.  At the end of this process, replace $\RR_{j+1}$ and $\QQ_{j+1}$ with $\widehat{\RR}_{j+1}$ and $\widehat{\QQ}_{j+1}$ to obtain the approximation
 \[ \|\bA_{j+1}\| \le \tol \cdot d_{j+1}\cdot\|\RR_{j+1}\|_F,\]
 since $d_{j+1} = d_j + 1$.

The error bound for the base case $j=1$ clearly holds, completing the induction.\sq

\medskip
The approximate QR factorization that results from this algorithm could be used directly for the approximation of leverage scores.  The perturbation theory of Ipsen and Wentworth~\cite{IW} describes how the tolerance in our algorithm will affect the accuracy of the resulting leverage scores.
We also note that for extra expediency this one-pass QR algorithm could be stopped when $\|\widehat{\RR}_j\|_F \approx \|\bA\|_F$ (at the cost of an extra pass through $\bA$ to compute $\|\bA\|_F$), or applied to only a random sampling of $k$ columns of $\bA$.  (Drma\v{c} and Gugercin propose a different random approach to DEIM index selection, based on sampling rows of $\bV$ to compute DEIM indices~\cite{DG}.)

\section{Computational Examples} \label{sec:examples}

This section presents some computational evidence illustrating the excellent approximation
properties of the DEIM-CUR factorization, consistent with the error analysis in Section~\ref{sec:theory}.  
For each of our three examples, we compare the accuracy of the DEIM-CUR factorization with several schemes based on leverage scores.  To remove random variations from our experiments, in most cases we  select columns and rows having the highest leverage scores; for the first example, we include results for random leverage score sampling.
For Example~1 we also study the effect of inaccurate singular vectors on the DEIM selection, and compare the accuracy of DEIM-CUR to CUR approximations based on the column-pivoted QR algorithm.

\newpage

\noindent \textbf{Example 1}.  \emph{Low-rank approximation of a sparse, nonnegative matrix}

\medskip
\noindent
The first example builds a matrix $\bA\in\IR^{300,000\times 300}$ of the form
\begin{equation} \label{eq:sparsex2}
\bA = \sum_{j=1}^{10} {2\over j}\,\bx_j^{}\by_j^T 
           + \sum_{j=11}^{300} {1\over j}\,\bx_j^{}\by_j^T,
\end{equation}
where $\bx_j \in \IR^{300,000}$ and $\by_j \in \IR^{300}$ are sparse vectors with random nonnegative entries (in MATLAB, $\bx_j = {\tt sprand(300000,1,0.025)}$ and $\by_j = {\tt sprand(300,1,0.025)}$).  In this instantiation, $\bA$ has 15,971,584 nonzeros, i.e., about 18\% of all entries are nonzero.  The form~(\ref{eq:sparsex2}) is not a singular value decomposition, since $\{\bx_j\}$ and $\{\by_j\}$ are not orthonormal sets; however, this decomposition suggests the structure of the SVD:  the singular values decay like $1/j$, and with the first ten singular values weighted more heavily to give a notable drop between $\sigma_{10}$ and $\sigma_{11}$.
We begin these experiments by computing $\bV$ and $\bW$ using 
MATLAB's economy-sized SVD routine (\verb|[V,S,W] = svd(A,'0')|).

\begin{figure}[b!]
\begin{center}
\includegraphics[scale=0.65]{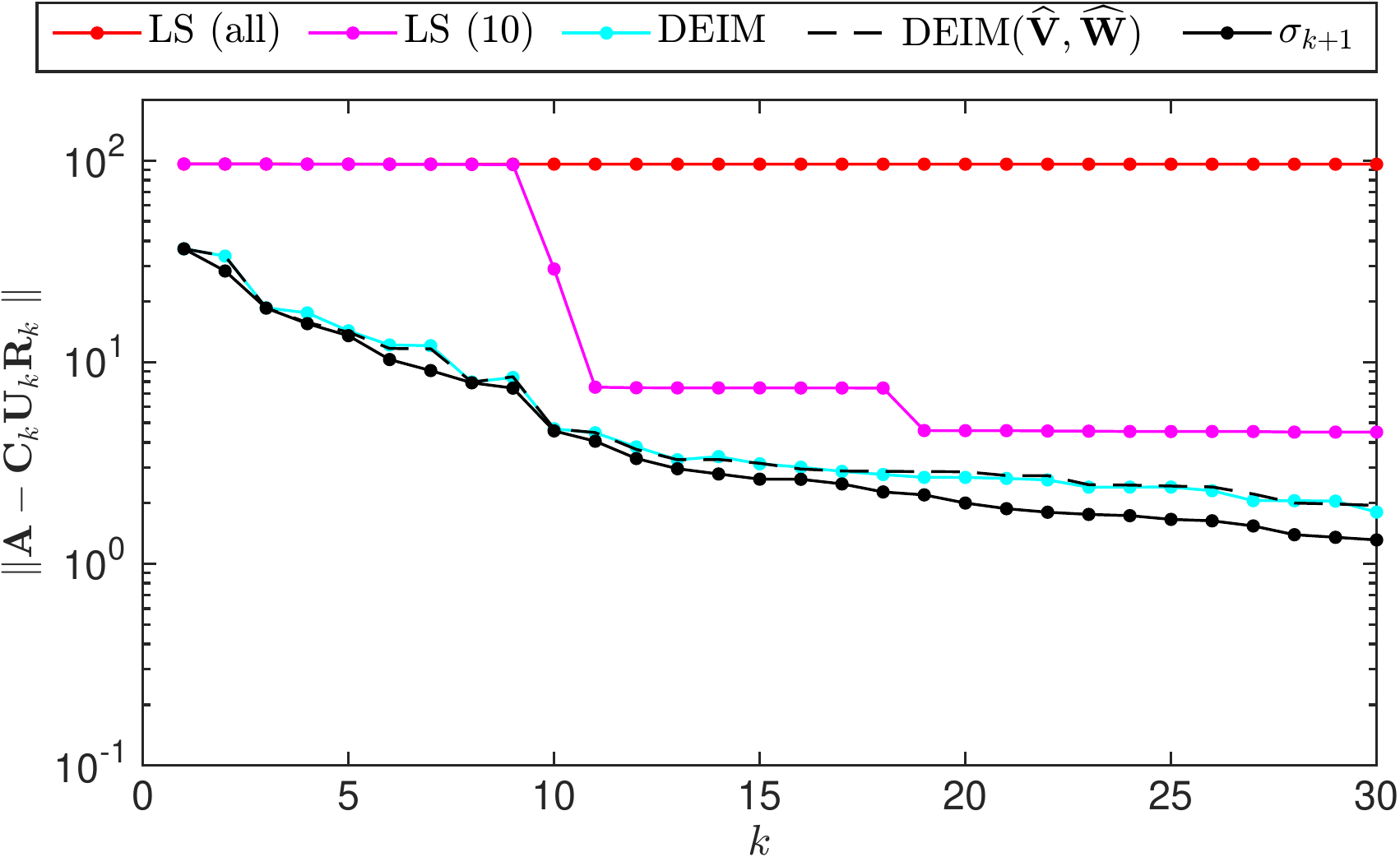}
\end{center}

\vspace*{-20pt}
\caption{\label{fig:sparse_ex4_cur}
Accuracy of CUR approximations for the sparse, nonnegative matrix~(\ref{eq:sparsex2}) using $k$ columns and rows, constructed by DEIM-CUR and two leverage score strategies:  ``LS (all)'' selects rows and columns with highest leverage scores computed using all 300~singular vectors; ``LS (10)'' only uses the leading ten singular vectors.  The ``DEIM($\widehat{\bV},\widehat{\bW}$)'' curve (nearly atop the ``DEIM'' curve) uses approximate singular vectors, described later.}
\end{figure}

Figure~\ref{fig:sparse_ex4_cur} compares the error $\|\bA-\bC\bU\bR\|$ for DEIM-CUR and methods that take $\bC$ and $\bR$ as the columns and rows of $\bA$ with the highest leverage scores.  These scores are computed using either all right and left singular vectors (300 of each), or using only the leading ten right and left singular vectors.  Both  approaches perform rather worse than DEIM-CUR, which closely tracks the optimal value $\sigma_{k+1}$.  

\begin{figure}[t!]
\begin{center}
\includegraphics[scale=0.45]{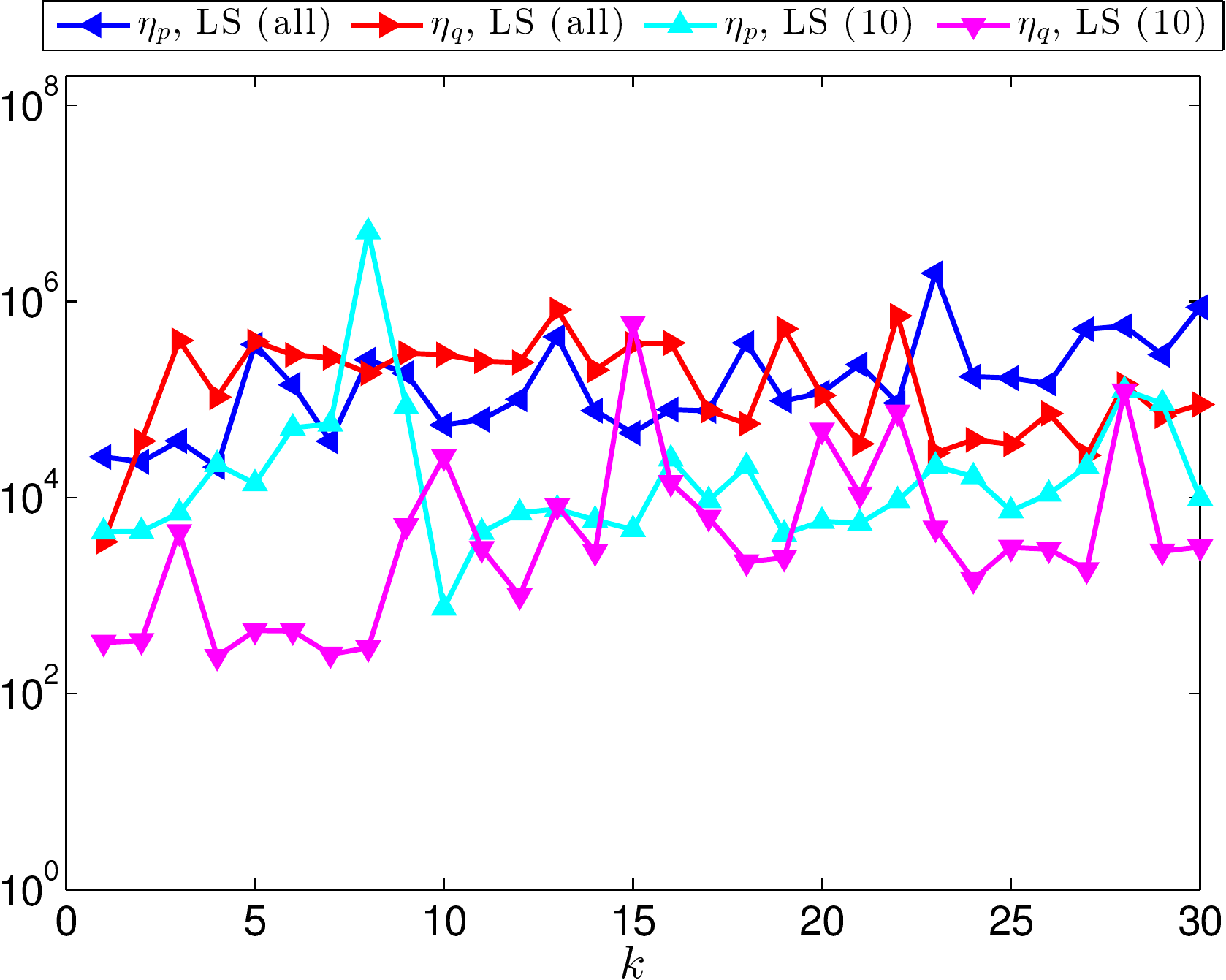}\quad
\includegraphics[scale=0.45]{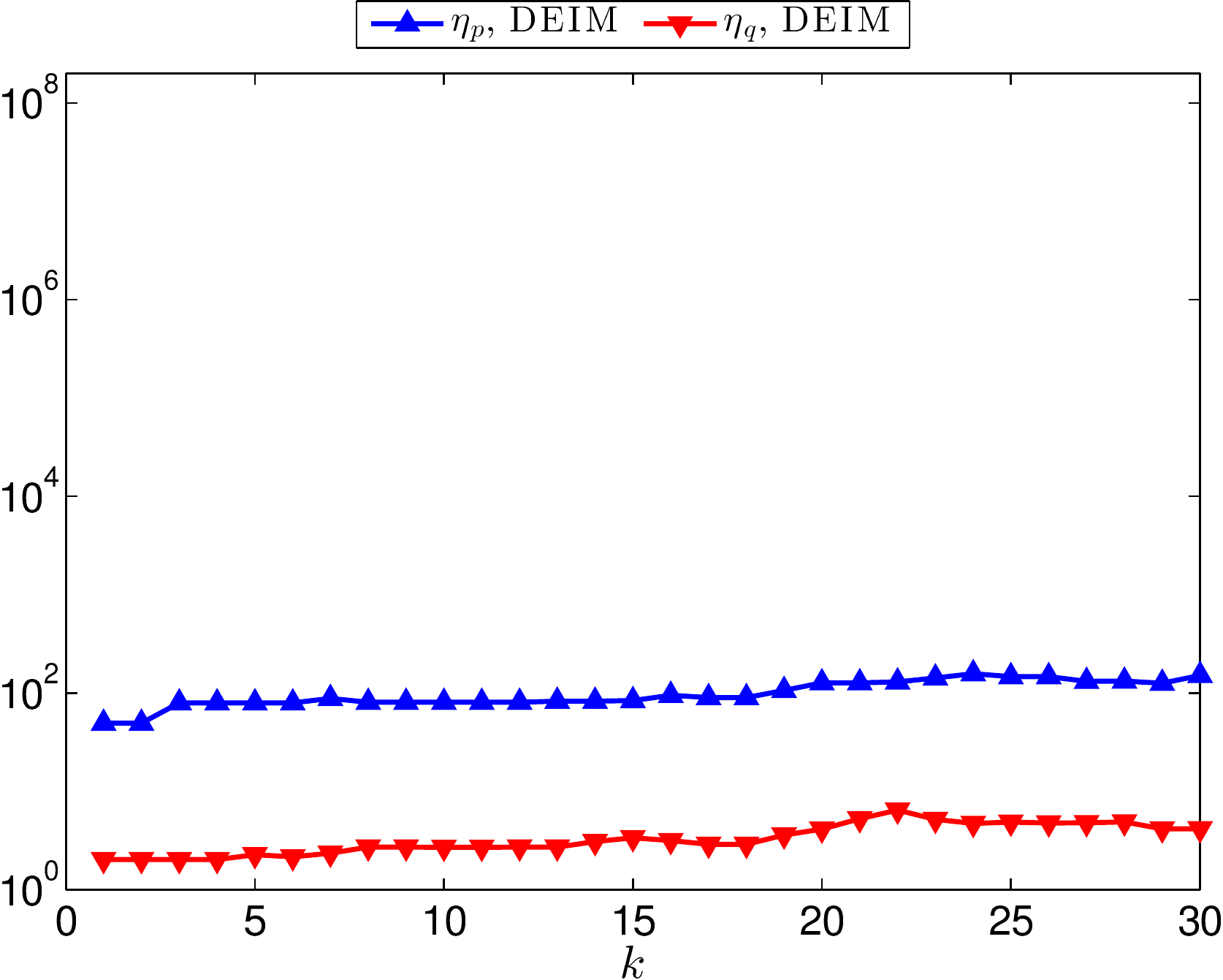}
\end{center}

\vspace*{-20pt}
\caption{\label{fig:sparse_ex4_eta}
Error constants $\eta_p = \|(\bP_k^T\bV_k)^{-1}\|$ and $\eta_q = \|(\bW_k^T\bQ)^{-1}\|$ for rows and columns selected using
 two leverage score strategies (left plot) and the DEIM algorithm (right plot), for the matrix~(\ref{eq:sparsex2}). }
\end{figure}

To gain insight into these results, we examine the interpolation constants $\eta_p$ and $\eta_q$ for all three approaches.  
Figure~\ref{fig:sparse_ex4_eta} shows that these  constants are largest for leverage scores based on all the singular vectors; using only ten singular vectors improves both the interpolation constants and the accuracy of the approximation (as seen in Figure~\ref{fig:sparse_ex4_cur}).  
The DEIM-CUR method gives better interpolation constants and more accurate approximations.  

\begin{figure}[t!]
\begin{center}
\includegraphics[scale=0.65]{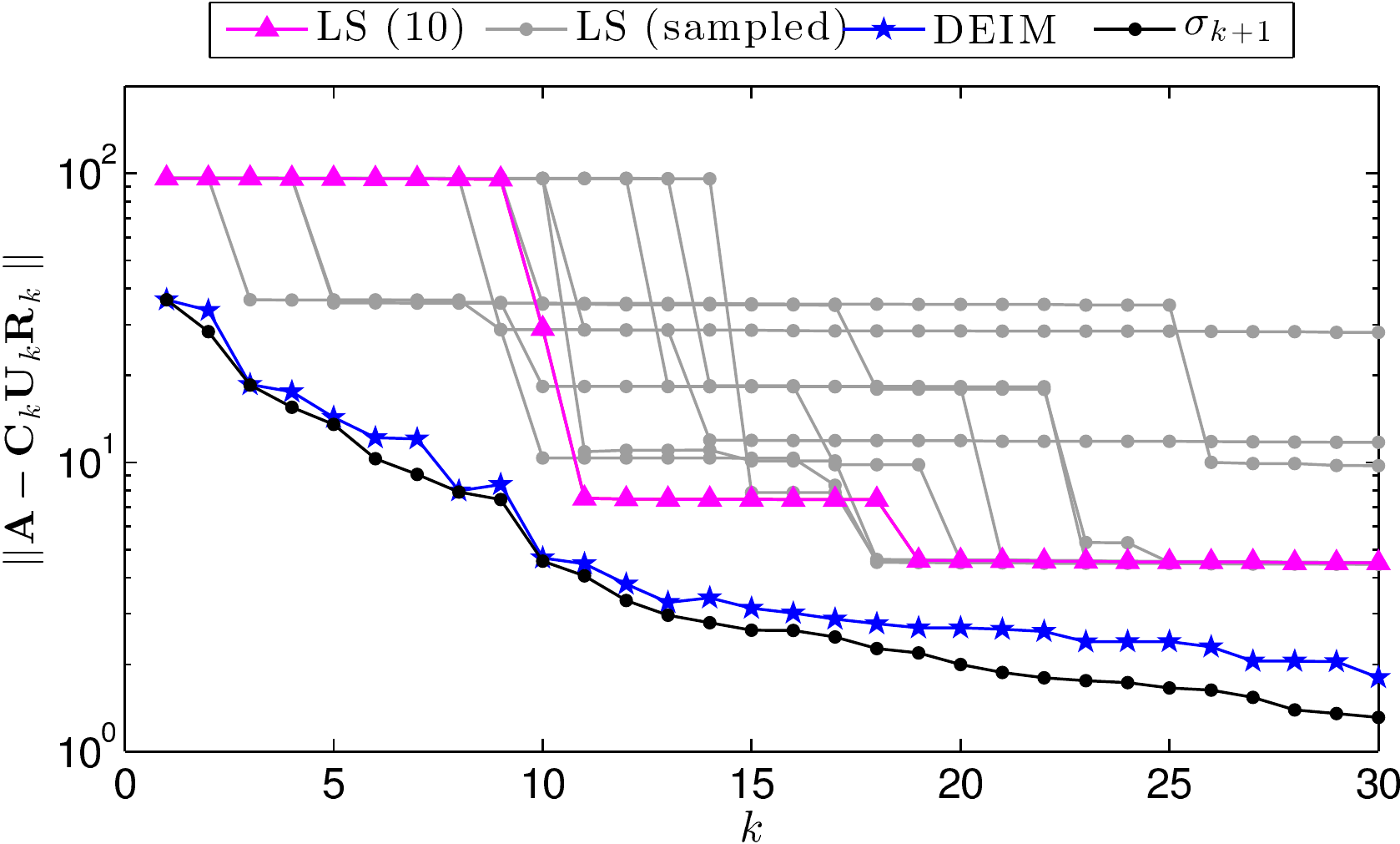}
\end{center}

\vspace*{-20pt}
\caption{\label{fig:sparse_ex5_cur}
Accuracy of CUR approximations for~(\ref{eq:sparsex2}) generated by randomly sampling rows and columns with probability weighted by leverage scores computed from the leading ten singular vectors.  All ten~trials (gray lines) perform similarly to the deterministic ``LS (10)'' approach, and worse than the DEIM-CUR approximation.}
\end{figure}

A CUR factorization can also be obtained by randomly sampling columns and rows of $\bA$, with the probability of selection weighted by leverage scores~\cite{MD09}.  We apply this approach on the current example, selecting $k=30$ rows and columns of $\bA$ with a probability given by the leverage scores computed from the leading ten singular vectors (normalized to give a probability distribution).  Figure~\ref{fig:sparse_ex5_cur} gives the results of ten independent experiments, showing that while sampling can sometimes yield better results than the deterministic leverage score approach, overall the approximations are still inferior to those from DEIM-CUR.

\begin{figure}[t!]
\begin{center}
\includegraphics[scale=0.45]{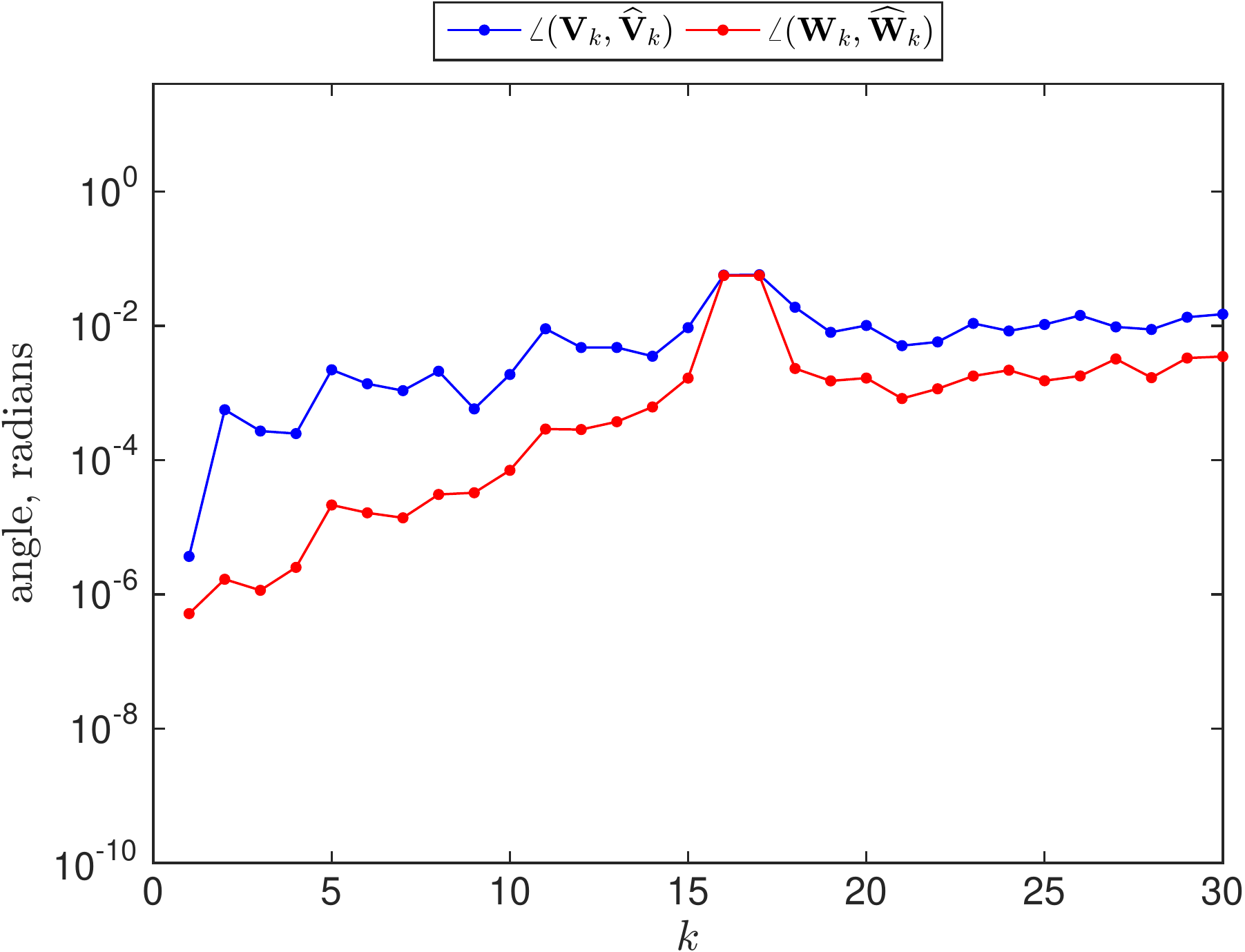}\quad
\includegraphics[scale=0.45]{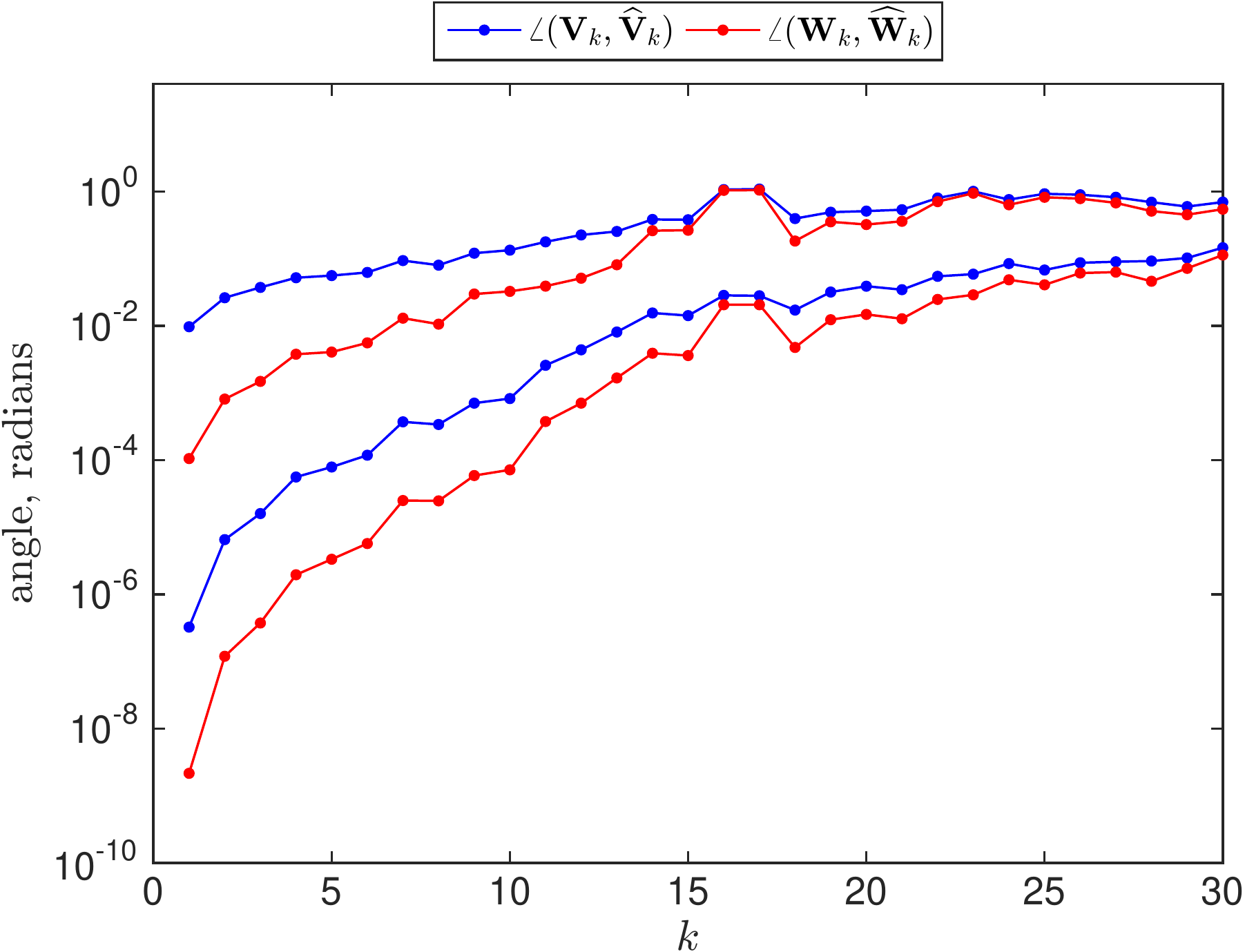}

\begin{picture}(0,0)
\put(-202,40){\small \emph{SVD via Incremental QR, ${\it tol}=10^{-4}$}}
\put(38,145){\rotatebox{10}{\small \emph{one application of $\bA$ and $\bA^{\kern-1pt T}$}}}
\put(47,60){\rotatebox{33}{\small \emph{two applications of $\bA$ and $\bA^{\kern-1pt T}$}}}
\put(45,40){\small \emph{Randomized SVD}}
\end{picture}
\end{center}

\vspace*{-32pt}
\caption{\label{fig:dirty_sparse4_angle}
The angle between the leading $k$-dimensional exact singular subspaces 
${\rm Ran}(\bV_k)$ and ${\rm Ran}(\bW_k)$
(generated by MATLAB's {\tt svd} command) and approximate
singular subspaces ${\rm Ran}(\widehat{\bV}_k)$ and ${\rm Ran}(\widehat{\bW}_k)$
for the matrix~(\ref{eq:sparsex2}).
On the left, $\widehat{\bV}_k$ and $\widehat{\bW}_k$ are generated using the Incremental QR algorithm described in Section~\ref{sec:qr}, with ${\it tol} = 10^{-4}$;
on the right, $\widehat{\bV}_k$ and $\widehat{\bW}_k$ are
generated using randomized SVD algorithm~\cite{HMT11}  using one and two applications of $\bA$ and $\bA^T$.}
\end{figure}
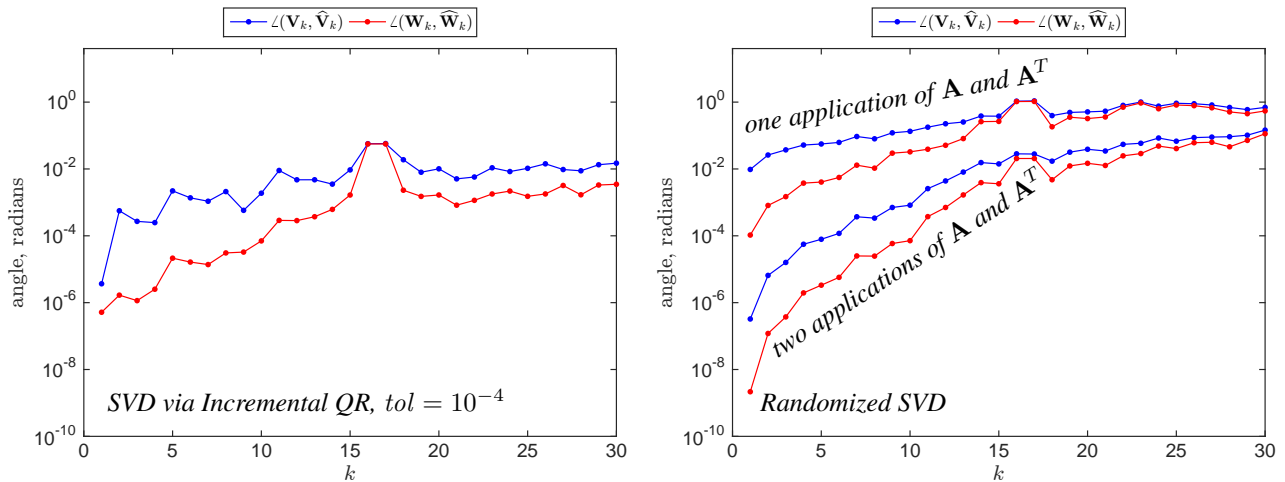

How robust is the DEIM-CUR approximation to errors in the singular vectors?
To investigate, we compute 
$\widehat{\bV}\approx \bV$ and $\widehat{\bW}\approx \bW$ using the
Incremental QR algorithm detailed in Section~\ref{sec:qr} (with ${\it tol}=10^{-4}$)
and the Randomized SVD algorithm described by
Halko, Martinsson, and Tropp~\cite[p.~227]{HMT11}.  
To give extreme examples of the latter, we compute $\widehat{\bV}$ and $\widehat{\bW}$ through 
\emph{only one or two} applications each of $\bA$ and $\bA^T$.%
\footnote{
This corresponds to $q=0$ and $q=1$ in the notation of~\cite[p.~227]{HMT11}.
Let the columns of $\bQ\in \IR^{m\times 2 k_{\rm max}}$ form an orthonormal basis for 
$(\bA\bA^T)^q\bA\bOmega$, where $\bOmega\in \IR^{n\times 2 k_{\rm max}}$ is a 
random matrix with i.i.d.\ Gaussian entries and we take $k_{\rm max} = 30$.
Then the leading $k_{\rm max}$~columns of $\bV$ and $\bW$ are approximated 
by taking the SVD of $\bQ^*\bA \in \IR^{2k_{\rm max} \times n}$. 
}
As Figure~\ref{fig:dirty_sparse4_angle} illustrates, in both cases the angle between the exact and approximate leading singular subspaces is significant, particularly as $k$ grows.
This drift in the subspaces has little effect on the accuracy of the DEIM approximations.

\begin{itemize}
\item The DEIM approximation using the Incremental QR algorithm is quite robust, choosing at most 3~different row indices and 2~different column indices for $k=1,\ldots,30$, with a relative discrepancy in $\|\bA-\bC_k\bU_k\bR_k\|$ of at most 9.27\% (and this realized only at step $k=30$).
\item When $\bA$ and $\bA^T$ are applied once in the Randomized SVD algorithm, the DEIM indices differ considerably from those drawn from exact singular vectors (e.g., for $k=30$, 20 of  30 row indices and 3 of 30 column indices differ), yet the quality of the approximation $\|\bA-\bC_k\bU_k\bR_k\|$ remains almost the same (relative difference of at most 10.45\%); see the dashed line in Figure~\ref{fig:sparse_ex4_cur}.
\item When $\bA$ and $\bA^T$ are applied twice, the DEIM indices are nearly identical (e.g., for $k=30$, 0 of 30 row indices and 2 of 30 column indices differ). On the scale of the plot in Figure~\ref{fig:sparse_ex4_cur}, $\|\bA-\bC_k\bU_k\bR_k\|$ could not be distinguished from the DEIM-CUR errors using exact singular vectors; the maximum relative discrepancy is 2.21\%.
\end{itemize}

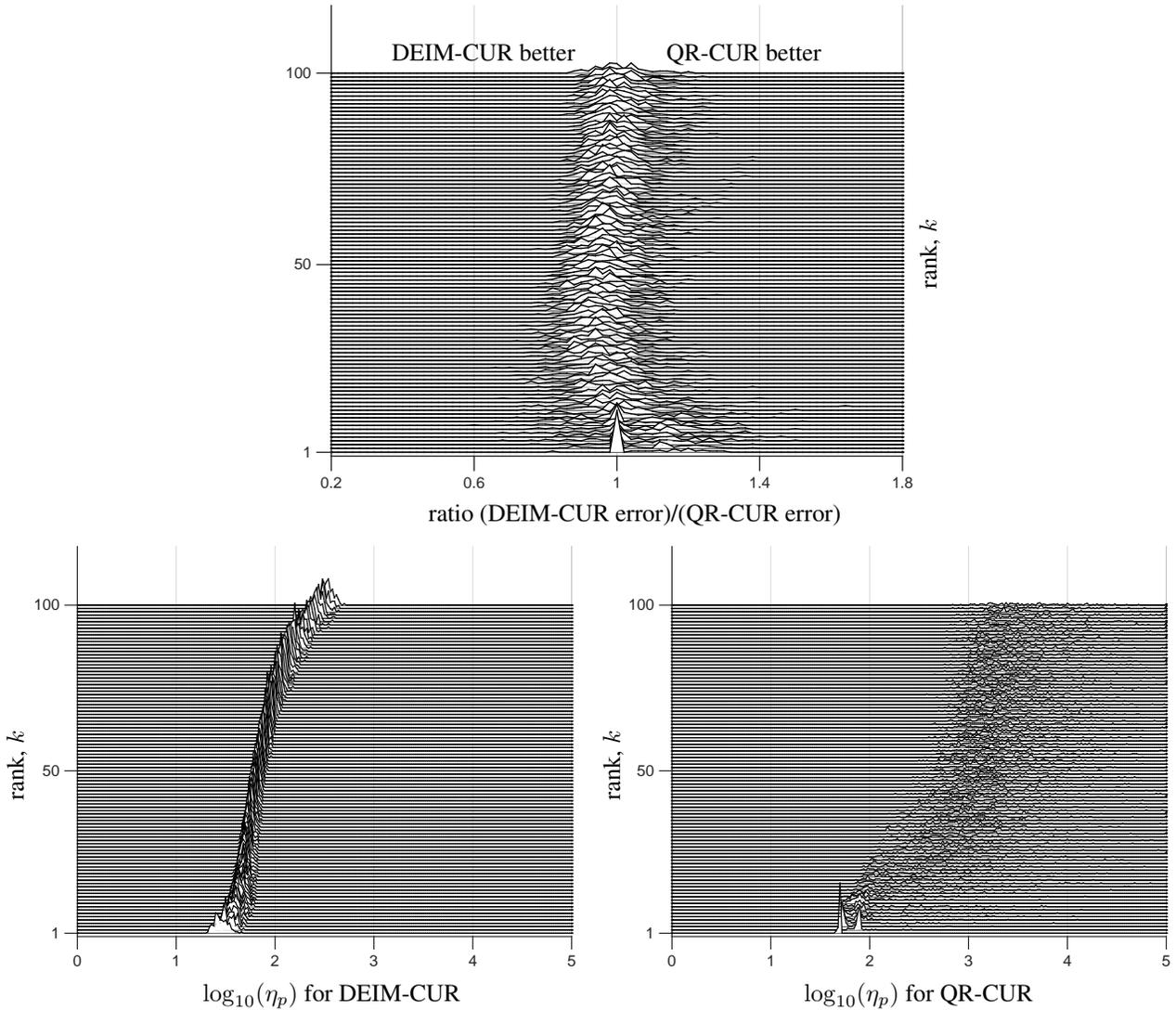
\begin{figure}[b!]
\begin{center}
\includegraphics[scale=0.52]{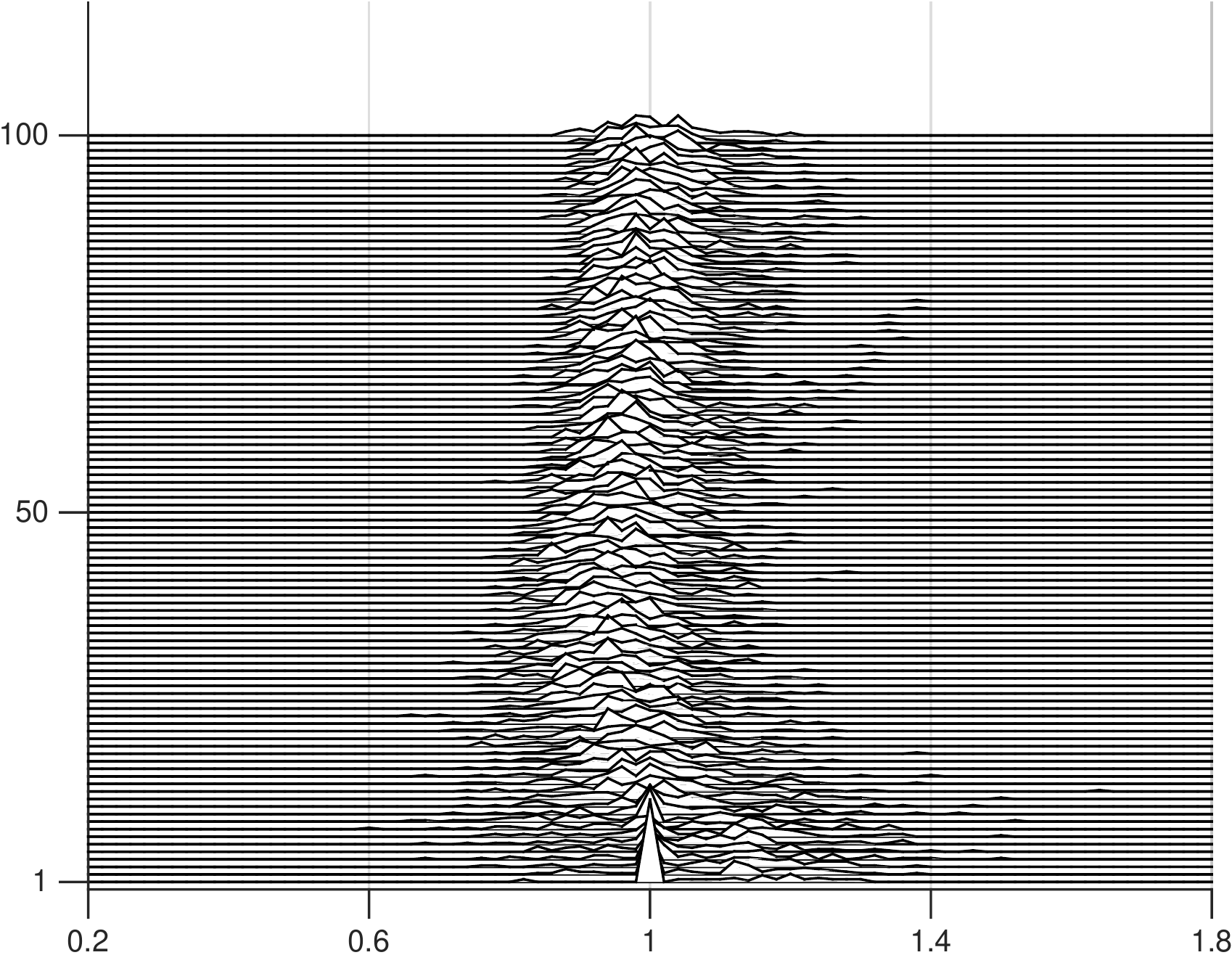}
\begin{picture}(0,0)
\put(-200,-13){\small ratio (DEIM-CUR error)/(QR-CUR error)}
\put(-215,175){\small DEIM-CUR better} 
\put(-103,175){\small QR-CUR better} 
\put(0,82){\small \rotatebox{90}{rank, $k$}}
\end{picture}

\vspace*{2em}
\includegraphics[scale=0.45]{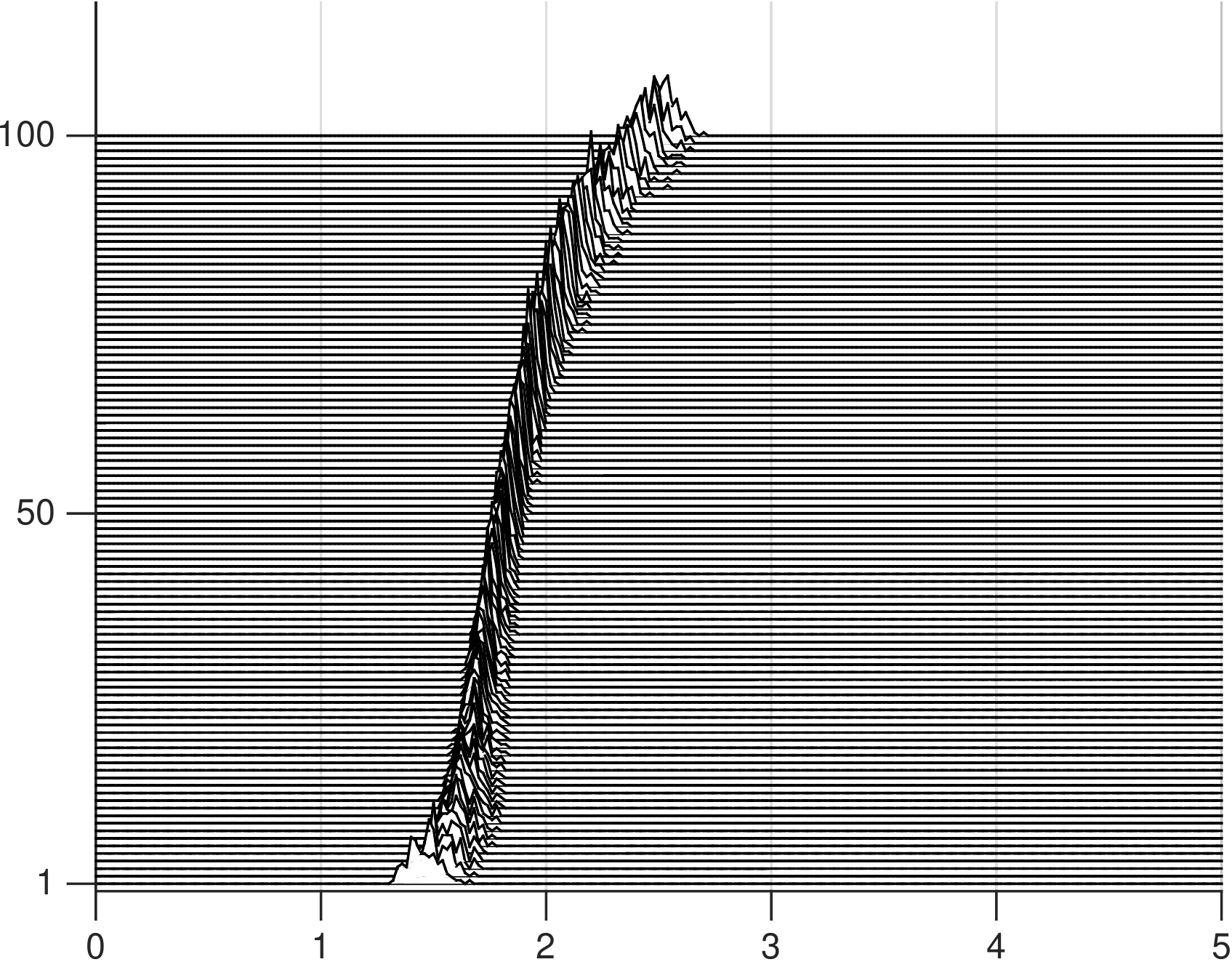}\qquad
\includegraphics[scale=0.45]{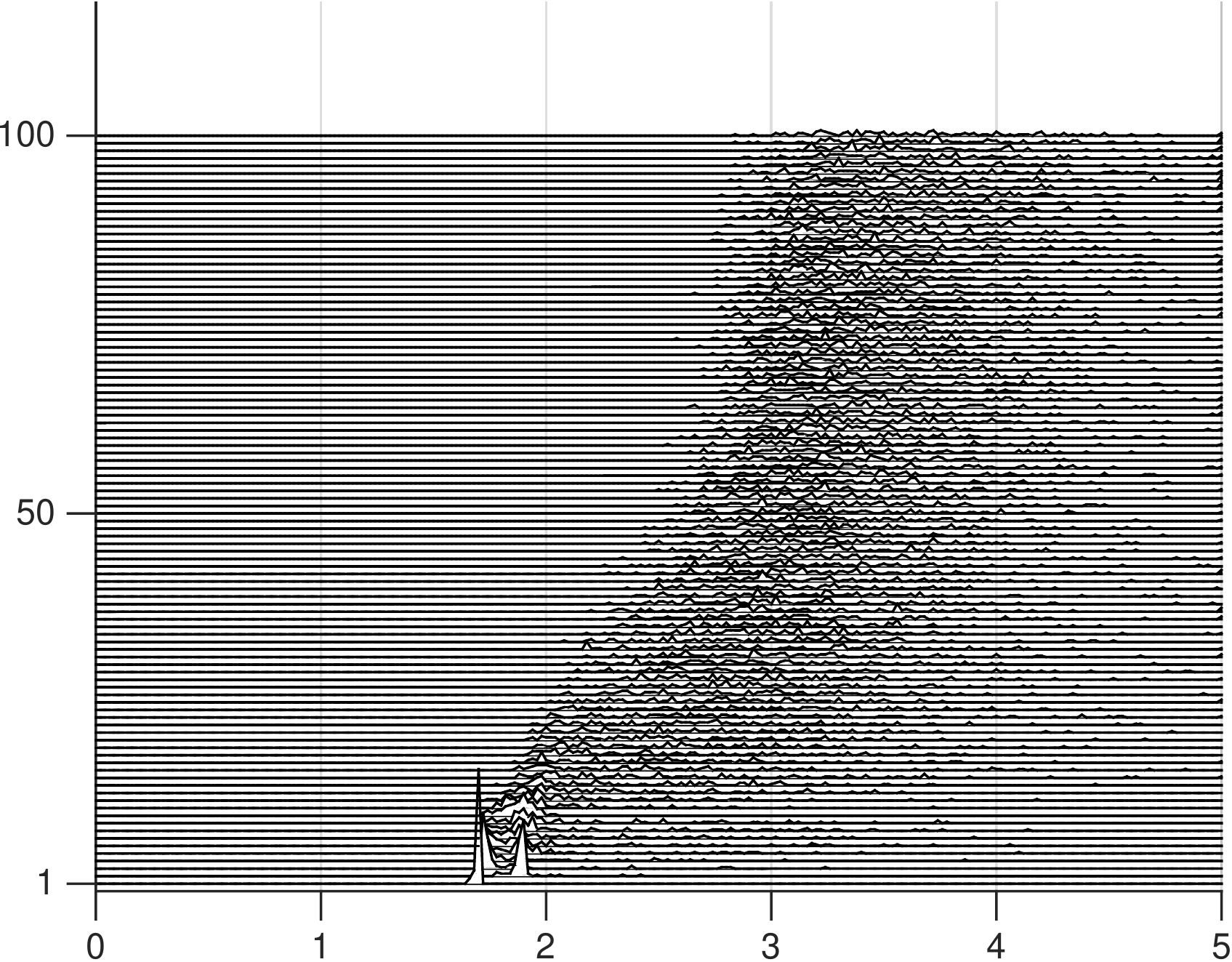}

\begin{picture}(0,0)
\put(-163,0){\small $\log_{10}(\eta_p)$ for DEIM-CUR}
\put( 82,0){\small $\log_{10}(\eta_p)$ for QR-CUR}
\put(-243,82){\small \rotatebox{90}{rank, $k$}}
\put(0,82){\small \rotatebox{90}{rank, $k$}}
\end{picture}
\end{center}

\vspace*{-10pt}
\caption{\label{fig:test10}
Comparison of DEIM-CUR and QR-CUR performance for
100~sparse random $300,000\times 300$ matrices of the form~(\ref{eq:sparsex2}).  
The top plot shows a histogram of the ratio of $\|\bA-\bC_k\bU_k\bR_k\|$ for DEIM-CUR and QR-CUR.  
The bottom plots compare the error constant $\eta_p = \|(\bP_k^T\bV_k)^{-1}\|$ for DEIM-CUR (left) and QR-CUR (right); 
note the logarithmic scale of the horizontal axes in the lower plots.} 
\end{figure}

Thus far we have only compared the DEIM-CUR approximations to CUR factorizations obtained from leverage scores, which also use singular vector information, thus illustrating how DEIM can use the same raw materials to better effect.
Next we compare DEIM-CUR to approximations computed using a different approach based on QR factorization of $\bA$; see, e.g., \cite{CGMR05,Ste99}.
Begin by computing a column-pivoted QR factorization of $\bA$; the first $k$ selected columns  give the indices $\bq$, from which we extract $\bC_k = \bA(@@:@@,\bq)$.
Next, a column-pivoted QR factorization of $\bC_k^T$ is performed; the first $k$ selected columns of $\bC_k^T$ 
give the indices $\bp$, from which we build $\bR_k = \bA(\bp,@@:@@)$.
We refer to this technique as ``QR-CUR.''

Figure~\ref{fig:test10} compares the results for 100~trials 
involving sparse random matrices of dimension $300,000\times 300$ 
having the form of our first experiment~(\ref{eq:sparsex2}).  DEIM-CUR and QR-CUR 
produce factorizations with similar accuracy, which we illustrate with a histogram of the ratio of $\|\bA-\bC_k\bU_k\bR_k\|$ for DEIM-CUR to QR-CUR, for $k=1,\ldots, 100$.  
(DEIM-CUR produces a smaller error when the ratio is less than one.)
While these errors are similar, the error constants $\eta_p$ and $\eta_q$ 
for the two methods are quite different.  
The bottom plots in Figure~\ref{fig:test10} compare histograms 
of $\log_{10} \eta_p$.  
For DEIM-CUR, the $\eta_p$ values are quite consistent
across the 100~random $\bA$, while for QR-CUR the $\eta_p$ values
are both larger and rather less consistent.  (The figures for $\eta_q$ are qualitatively identical, but about an order of magnitude smaller for both methods.)  

\begin{figure}[ht!]
\begin{center}
\includegraphics[scale=0.65]{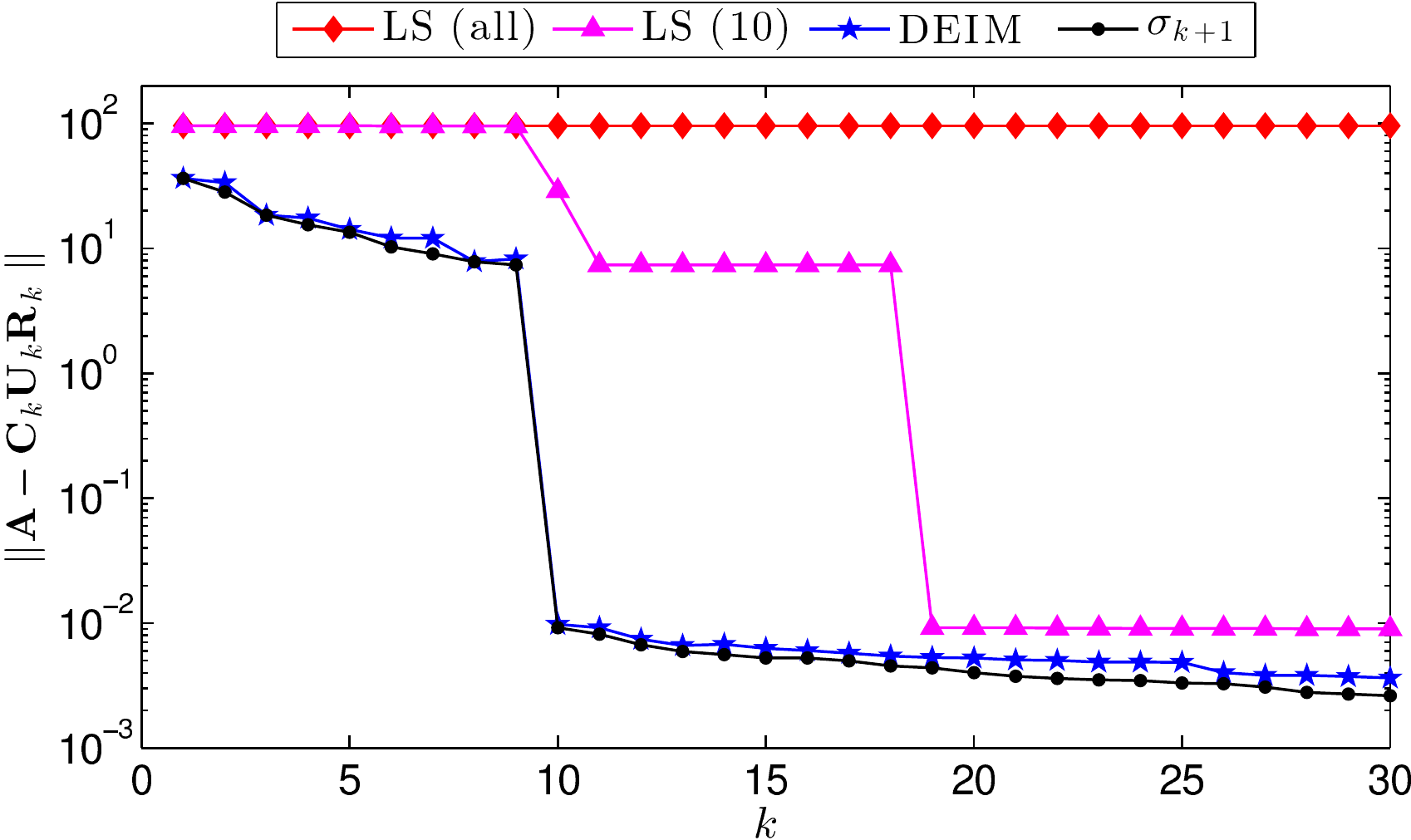}
\end{center}

\vspace*{-20pt}
\caption{\label{fig:sparse_ex3_cur}
Accuracy of CUR approximations using $k$ rows and columns, for DEIM-CUR and two leverage score strategies for the sparse, nonnegative matrix~(\ref{eq:sparsex1000}).  ``LS (all)'' selects rows and columns having the highest leverage scores computed using all 300~singular vectors; ``LS (10)'' uses the leading 10~singular vectors.}
\end{figure}

\begin{figure}[ht!]
\begin{center}
\includegraphics[scale=0.45]{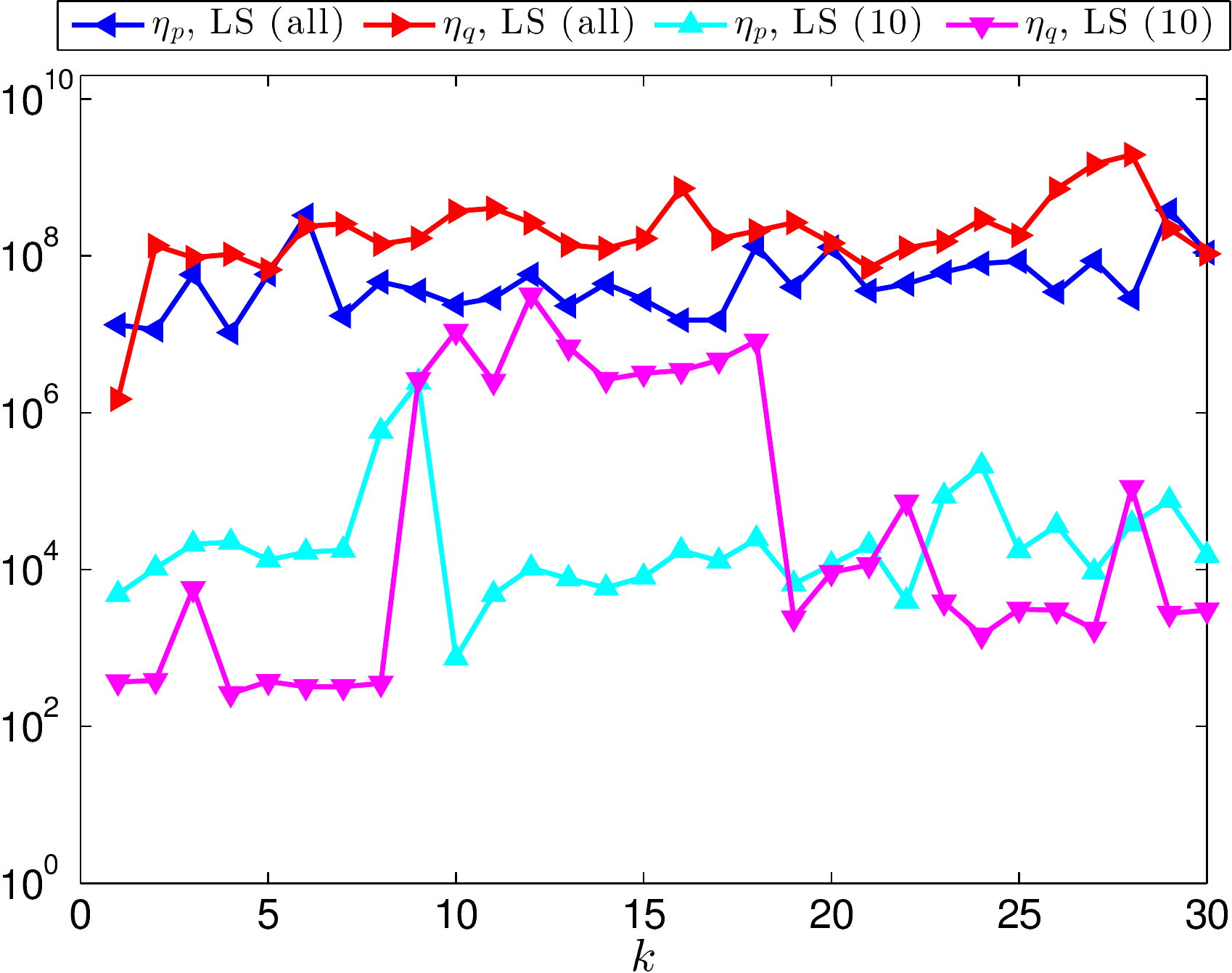}\quad
\includegraphics[scale=0.45]{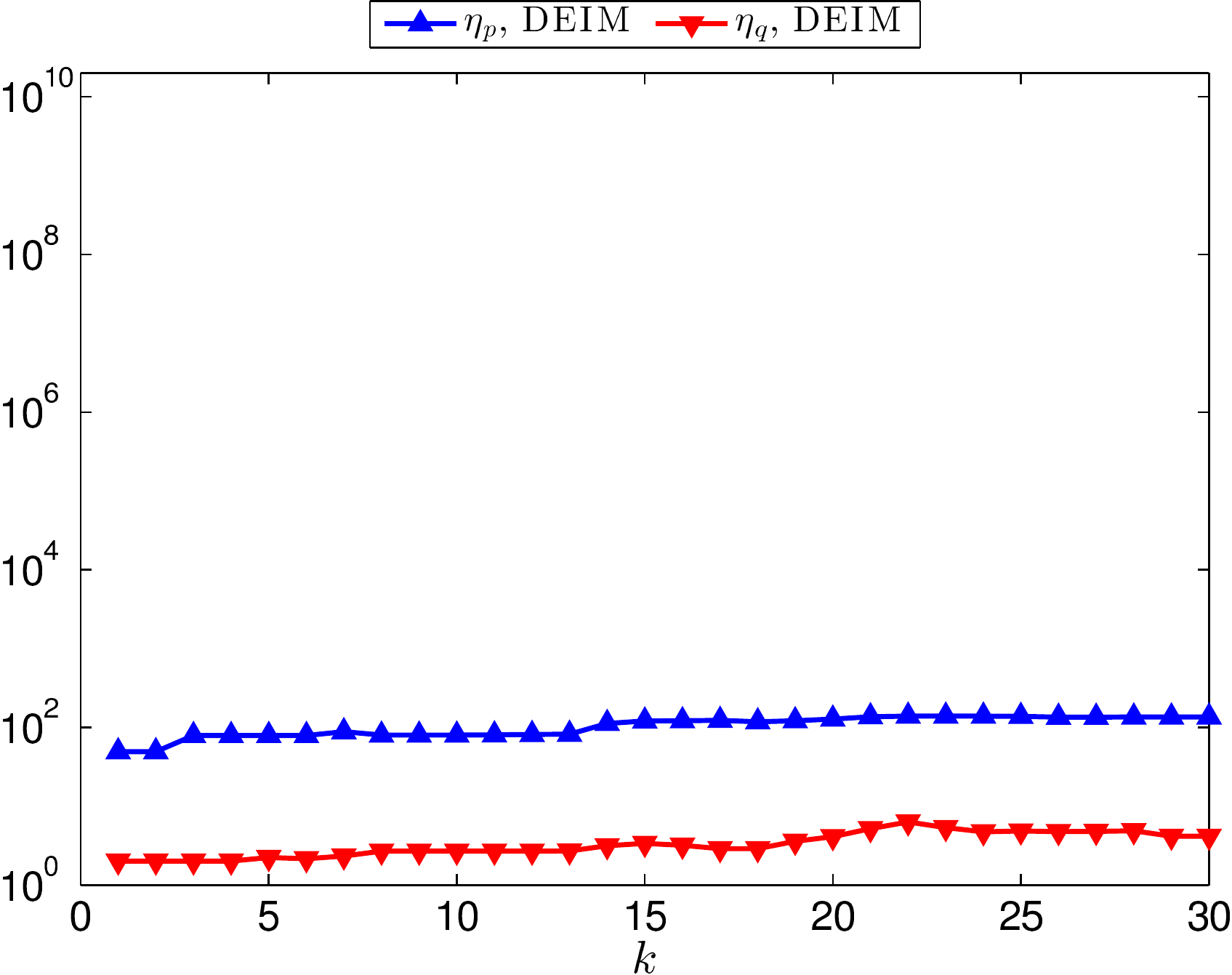}
\end{center}

\vspace*{-20pt}
\caption{\label{fig:sparse_ex3_eta}
Error constants $\eta_p = \|(\bP_k^T\bV_k)^{-1}\|$ and $\eta_q = \|(\bW_k^T\bQ_k)^{-1}\|$ for rows and columns selected using two 
leverage score strategies (left plot) and the DEIM algorithm (right plot), for the sparse matrix $\bA$ given in~(\ref{eq:sparsex1000}). }
\end{figure}

\medskip
The advantage of DEIM-CUR over approximations based on leverage scores remains when the singular values decrease more sharply.
Modify~(\ref{eq:sparsex2}) to give a more significant drop between $\sigma_{10}$ and $\sigma_{11}$:
\begin{equation} \label{eq:sparsex1000}
\bA = \sum_{j=1}^{10} {1000\over j}\,\bx_j^{}\by_j^T 
           + \sum_{j=11}^{300} {1\over j}\,\bx_j^{}\by_j^T.
\end{equation}
As seen in Figures~\ref{fig:sparse_ex3_cur} and~\ref{fig:sparse_ex3_eta}, the DEIM-CUR approach again delivers excellent approximations, while selecting the rows and columns with highest leverage scores does not perform nearly as well.  (In Figure~\ref{fig:sparse_ex3_eta}, note the significant jump in the ``LS~(10)'' error constant $\eta_q$ corresponding to those $k$ values where $\|\bA-\bC_k\bU_k\bR_k\|/\sigma_{k+1}$ is large.)
 
\newpage
\noindent \textbf{Example 2}.  \emph{TechTC term document data}

\medskip
\noindent
The second example, adapted from Mahoney and Drineas~\cite{MD09}, computes the CUR factorization of a term document matrix with data drawn from the Technion Repository of Text Categorization Datasets (TechTC)~\cite{GM04}.  The rows of the data matrix correspond to websites (consolidated from multiple webpages), while the columns correspond to ``features'' (words from the text of the webpages).  The $(j,k)$ entry of $\bA$ reflects the importance of the feature text on the given website; most entries are zero. For this experiment we use TechTC-100 test set~26, which concatenates a data set relating to Evansville, Indiana (id~10567) with another for Miami, Florida (id~11346).  Following Mahoney and Drineas~\cite{MD09}, we omit all features with four or fewer characters from the data set, leaving a matrix with 139~rows and 15,170~columns.  Each row of $\bA$ is then scaled to have unit 2-norm.  Ideally a CUR factorization not only gives an accurate low-rank approximation to $\bA$, but also selects rows corresponding to representative webpages from each geographic area, and columns corresponding to meaningful features.  

\begin{figure}[t!]
\begin{center}
\includegraphics[scale=0.65]{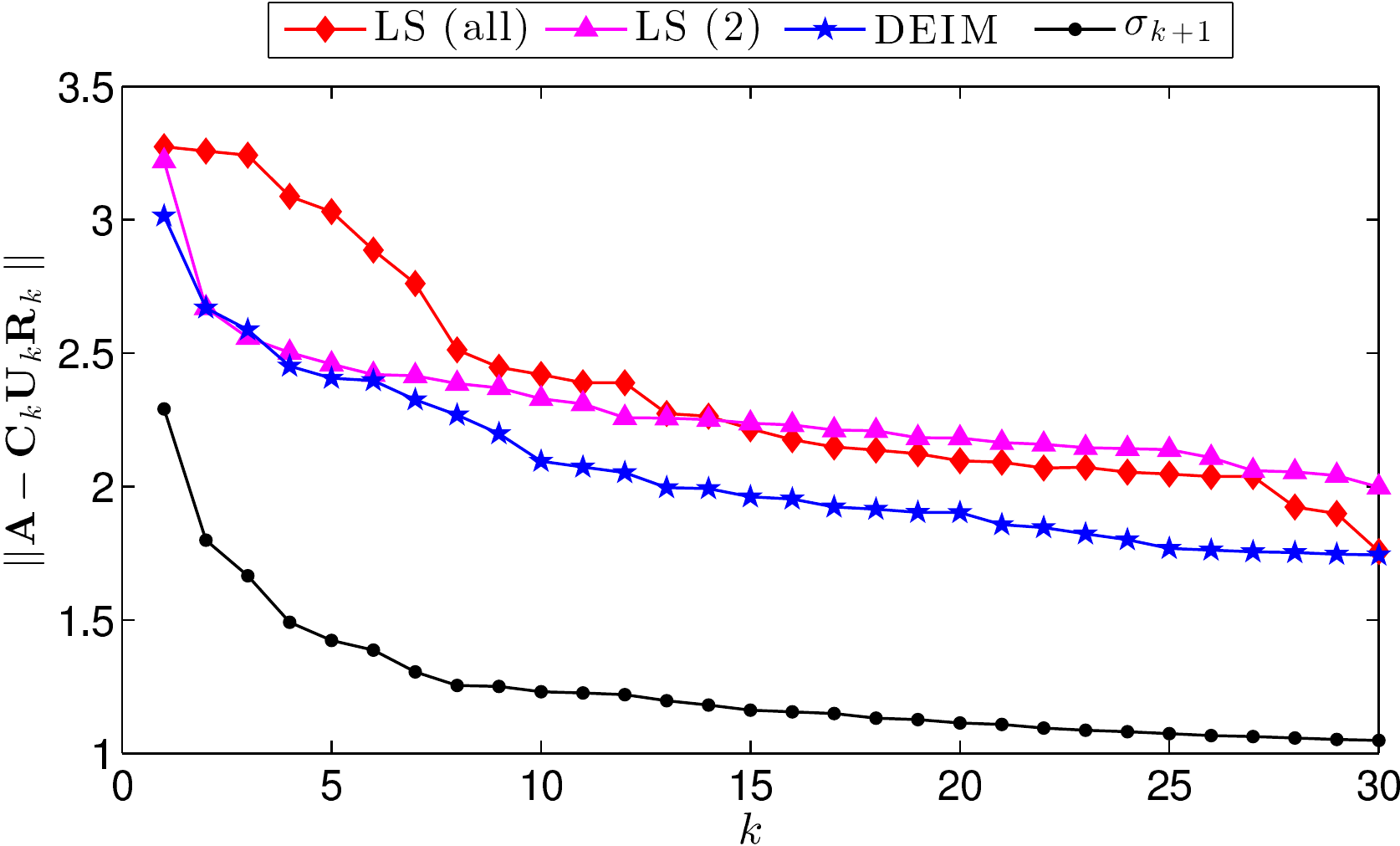}
\end{center}

\vspace*{-20pt}
\caption{\label{fig:techtc_rs_cur} Accuracy of CUR factorizations for the TechTC example, selecting rows and columns using top leverage scores for all singular vectors and the leading two singular vectors, and DEIM.}
\end{figure}

\begin{figure}[t!]
\begin{center}
\includegraphics[scale=0.65]{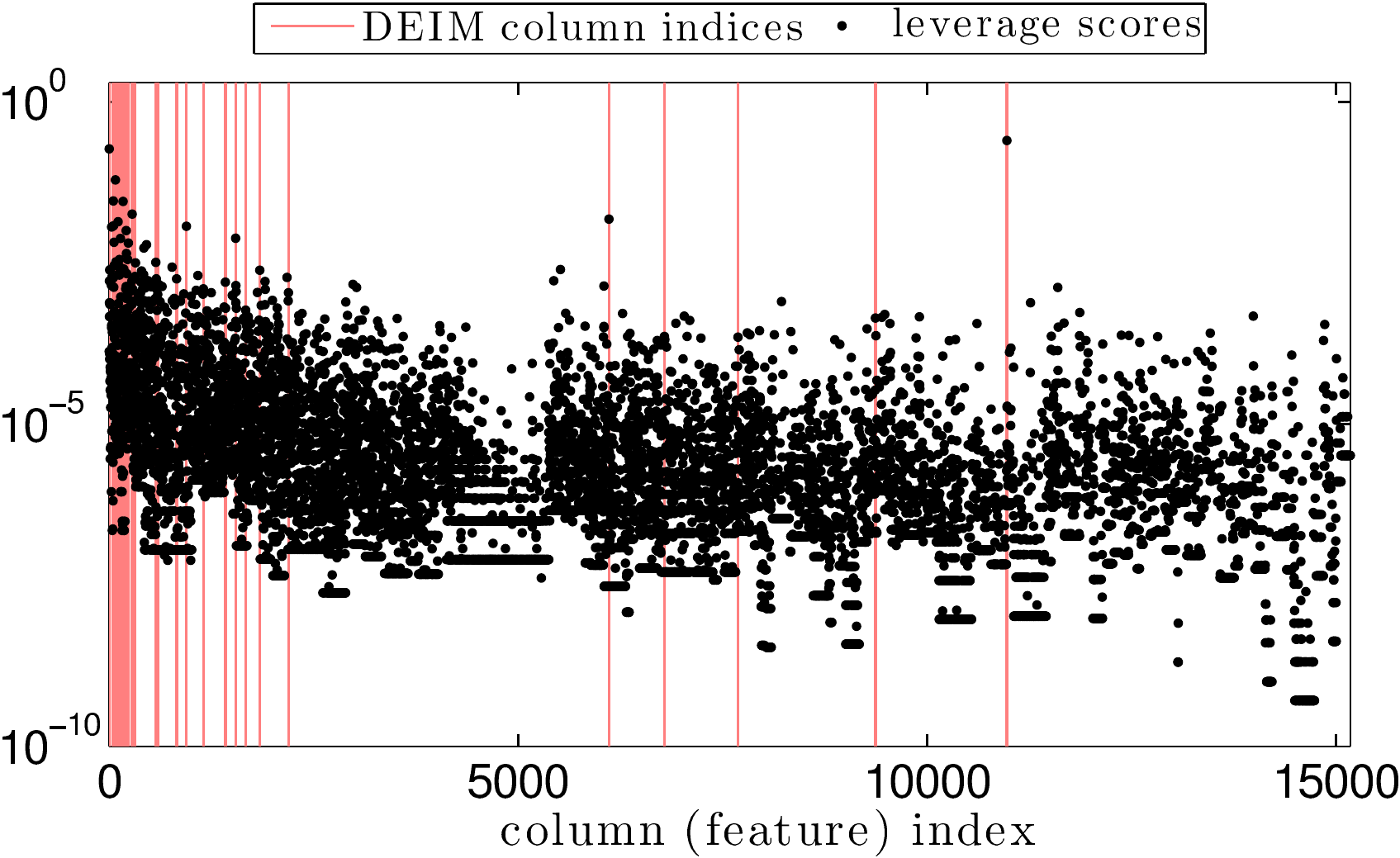}
\end{center}

\vspace*{-20pt}
\caption{\label{fig:techtc_rs_comp_ls2} The columns selected by DEIM for the TechTC example, compared to leverage scores from the leading two singular vectors.}
\vspace*{-10pt}
\end{figure}

Figure~\ref{fig:techtc_rs_cur} compares DEIM-CUR approximations to row and column selection based on highest leverage scores (from all singular vectors, or the two leading singular vectors).  The DEIM-CUR approximations are typically more accurate than those based on leverage scores, but all approaches give errors roughly two times larger than the slowly-decaying optimal value of $\sigma_{k+1}$. 
How do the DEIM columns (features) compare to those with the highest leverage scores?  Figure~\ref{fig:techtc_rs_comp_ls2} shows the leverage scores associated with each column of $\bA$ (based on the two leading singular vectors), along with the first 30~columns selected by DEIM.  While the columns with highest leverage scores were found by DEIM, there are DEIM columns with marginal leverage scores, and vice versa.  This data is more easily parsed in Table~\ref{tbl:techtc}, which lists the features corresponding to the first 20~DEIM columns.  (To ease comparison, we normalize leverage scores so that the maximum value is one.)  The leading features identified by DEIM, including ``evansville'' (first DEIM point), ``florida'' (second), ``miami'' (sixth), and ``indiana'' (nineteenth), indeed reveal key geographic terms.  These terms scored at least as high when ranked by leverage scores based on two leading singular vectors; when all singular vectors are used, the scores of these terms generally drop, relative to other features.  Overall, one notes that DEIM selects a significantly different set of indices than those valued by leverage scores, and, as seen in Figure~\ref{fig:techtc_rs_cur}, tends to provide a somewhat better low-rank approximation.

\begin{table}[t!]
\caption{\label{tbl:techtc}
The features selected by DEIM-CUR for the TechTC data set, compared to the (scaled) leverage scores using the leading two singular vectors, and all singular vectors.}

\vspace*{-5pt}
\begin{center}\small
\begin{tabular}{cr|cc|cc|l} \hline
\multicolumn{1}{c}{DEIM } &
\multicolumn{1}{c}{index} &
\multicolumn{2}{|c}{LS (2)} &
\multicolumn{2}{|c}{LS (all) } &
\multicolumn{1}{|c}{} \\
\multicolumn{1}{c}{rank } &
\multicolumn{1}{c}{$q_j$} &
\multicolumn{1}{|c}{rank } &
\multicolumn{1}{c}{score} &
\multicolumn{1}{|c}{rank } &
\multicolumn{1}{c}{score} &
\multicolumn{1}{|c}{feature} \\ \hline
     1 &      10973 &      1 &  1.000 &      4 &  0.875 & evansville\\
     2 &          1 &      2 &  0.741 &      8 &  0.726 & florida\\
     3 &       1547 &     13 &  0.031 &      2 &  0.948 & spacer\\
     4 &        109 &      8 &  0.055 &     66 &  0.347 & contact\\
     5 &        209 &     12 &  0.040 &     32 &  0.458 & service\\
     6 &         50 &      4 &  0.116 &      6 &  0.739 & miami\\
     7 &        824 &     46 &  0.007 &      5 &  0.809 & chapter\\
     8 &       1841 &     33 &  0.010 &     20 &  0.537 & health\\
     9 &        171 &      5 &  0.113 &     13 &  0.617 & information\\
    10 &        234 &     16 &  0.026 &     37 &  0.436 & events\\
    11 &        595 &     84 &  0.004 &     15 &  0.576 & church\\
    12 &         60 &     15 &  0.026 &     67 &  0.347 & email\\
    13 &        945 &     10 &  0.047 &     30 &  0.474 & services\\
    14 &       1670 &    129 &  0.002 &      1 &  1.000 & bullet\\
    15 &        216 &     35 &  0.009 &     38 &  0.430 & music\\
    16 &         78 &      3 &  0.246 &     24 &  0.492 & south\\
    17 &        213 &     19 &  0.018 &    110 &  0.259 & their\\
    18 &        138 &     14 &  0.030 &     43 &  0.408 & please\\
    19 &       6110 &      7 &  0.060 &     95 &  0.280 & indiana\\
    20 &       1152 &     70 &  0.005 &    152 &  0.221 & member
\end{tabular}\end{center}
\vspace*{-2em}
\end{table}

\medskip
\noindent \textbf{Example 3}.  \emph{Tumor detection in genetics data}

\medskip
\noindent
Our final example uses the GSE10072 cancer genetics data set from the National Institutes of Health, previously investigated by Kundu, Nambirijan, and Drineas~\cite{KND}.  The matrix $\bA\in\IR^{22,283\times 107}$ contains data for 22,283~probes applied to 107~patients.  The $(j,k)$ entry of $\bA$ reflects how strongly patient $k$ responded to probe $j$.  This experiment seeks probes that segment the population into two clusters: the 58~patients with tumors, and the 49~without.%
\footnote{The data is available from \url{http://www.ncbi.nlm.nih.gov/geo/query/acc.cgi?acc=GSE10072}.}
To center the data, we subtract the mean of each row from all entries in that row.  
As shown in~\cite{KND}, the leading two principal vectors of this matrix segment the population very well. 

\begin{figure}[ht]
\begin{center}
 \includegraphics[scale=0.65]{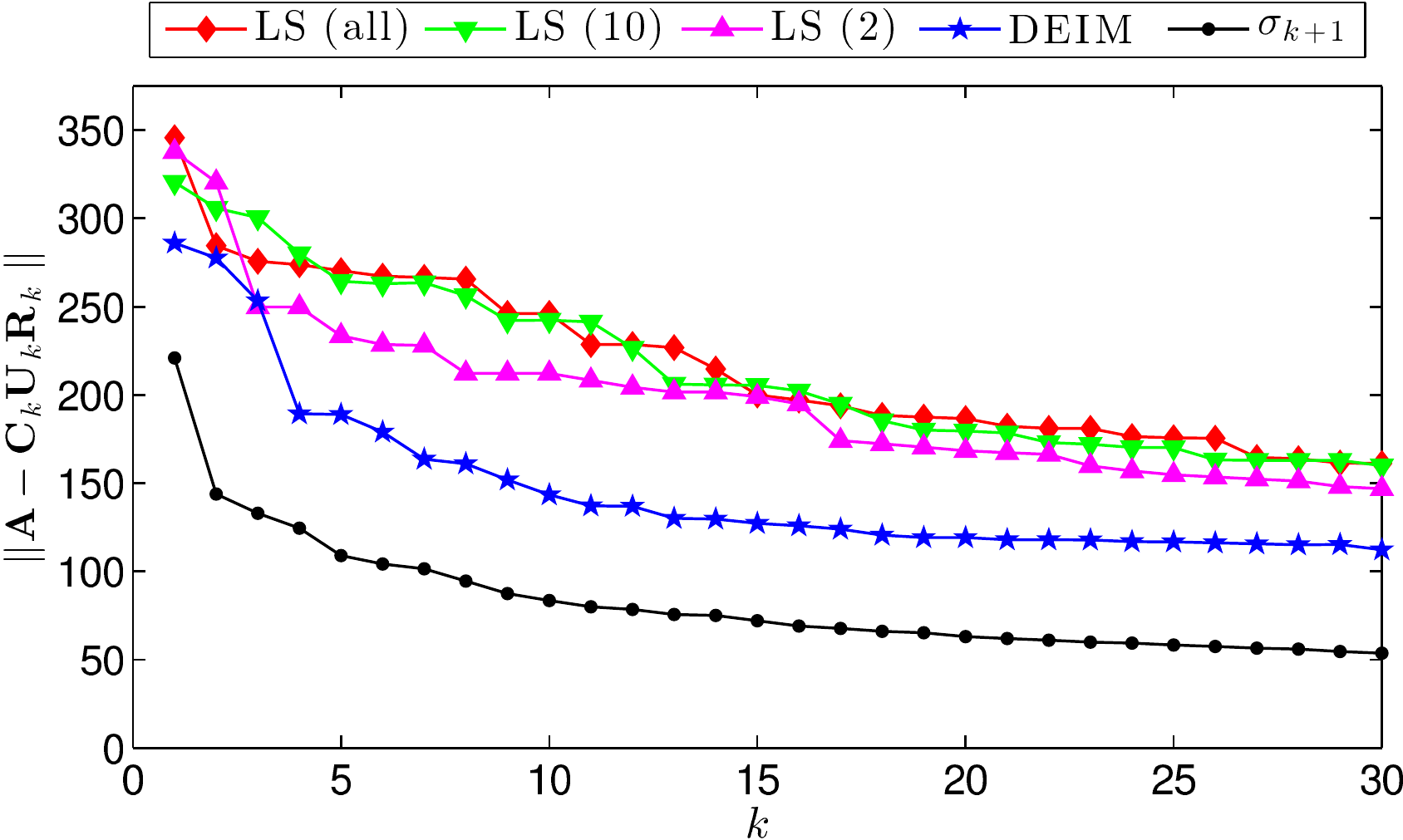}
 \end{center}
 
 \vspace*{-17pt}
 \caption{\label{fig:gene_cur}
 Accuracy of CUR factorizations for a genetics data set.
DEIM-CUR consistently outperforms factorizations
 derived by taking the rows and columns with largest leverage scores,
 regardless of whether these scores are drawn from all singular vectors,
 the leading ten singular vectors, or the leading two singular vectors.
 }
 \end{figure}

Like the TechTC data, the singular values of $\bA$ decay slowly, as seen in Figure~\ref{fig:gene_cur}.  Once again the DEIM-CUR procedure produces a more accurate low-rank approximation than obtained by selecting the rows and columns with highest leverage scores, whether those are computed using all the singular vectors, or just the leading two or ten.


Table~\ref{tbl:gene_deim} reports the first 15~rows selected by the DEIM-CUR process, along with the corresponding leverage scores based on two, ten, and all singular vectors.  Do the probes selected by DEIM  discriminate the patients with tumors (``sick'') from those without (``well'')? 
To investigate, for each selected probe we count the number of large positive entries corresponding to sick and well patients.%
\footnote{In particular, we call an entry of the mean-centered matrix $\bA$ large if its value exceeds one.  Of the 22,283 probes, for only 23~probes do at least~30 of the~58 sick patients have such large entries; for only~95 probes do at least~30 of the~49 well patients have large entries.  There is no overlap between the probes that are strongly expressed by the sick and well patients.}  Some but not all of the DEIM-CUR probes effectively select only sick or well patients.
Contrast these results with Table~\ref{tbl:gene_ls2}, which shows the probes with highest leverage scores (drawn from the leading two singular vectors).  Only four of these probes were also selected by the DEIM procedure (even if we continue the DEIM procedure to select the maximum number, $n=107$, of indices).  This discrepancy is quite different from the good agreement between DEIM and leverage score indices for the TechTC data in Table~\ref{tbl:techtc},  despite the similar dimensions and the comparably slow decay of the singular values.

\begin{table}[t!]
\caption{\label{tbl:gene_deim}
Genetics example: the probes selected by DEIM-CUR, compared to the (scaled) leverage scores using the leading two singular vectors, ten singular vectors, and all singular vectors.}
\begin{center}\small
\begin{tabular}{cr|ll|cc|cc|cc|cc} \hline
\multicolumn{1}{c}{DEIM } &
\multicolumn{1}{c}{index} &
\multicolumn{1}{|c}{probe} &
\multicolumn{1}{c}{gene} &
\multicolumn{1}{|c}{number} &
\multicolumn{1}{c}{number} &
\multicolumn{2}{|c}{LS (2) } &
\multicolumn{2}{|c}{LS (10)} &
\multicolumn{2}{|c}{LS (all) } \\
\multicolumn{1}{c}{rank } &
\multicolumn{1}{c}{$q_j$} &
\multicolumn{1}{|c}{set } &
\multicolumn{1}{c}{name} &
\multicolumn{1}{|c}{sick } &
\multicolumn{1}{c}{well} &
\multicolumn{1}{|c}{rank } &
\multicolumn{1}{c}{score} &
\multicolumn{1}{|c}{rank } &
\multicolumn{1}{c}{score} &
\multicolumn{1}{|c}{rank } &
\multicolumn{1}{c}{score}  \\ \hline
     1 &       9565 & 210081\_at &  AGER &   2 &  45 &     1 &  1.000 &    45 &  0.504 &    386 &  0.123 \\
     2 &      14270 & 214895\_s\_at & ADAM10  &   8 &   3 &   211 &  0.173 &  1171 &  0.108 &   3344 &  0.036 \\
     3 &       8650 & 209156\_s\_at & COL6A2  &   5 &   6 & 15156 &  0.005 &   252 &  0.245 &    708 &  0.091 \\
     4 &      11057 & 211653\_x\_at & AKR1C2  &  18 &   1 &  6440 &  0.017 &    11 &  0.656 &    146 &  0.185 \\
     5 &      14153 & 214777\_at & IGKV4-1  &  27 &   3 &   281 &  0.148 &    19 &  0.607 &    106 &  0.209 \\
     6 &      18976 & 219612\_s\_at & FGG  &  17 &  17 &  2591 &  0.039 &     2 &  0.956 &      4 &  0.825 \\
     7 &       3831 & 204304\_s\_at & PROM1  &  16 &   4 &   992 &  0.073 &    70 &  0.417 &     32 &  0.345 \\
     8 &       3351 & 203824\_at & TSPAN8  &  17 &   4 &  9687 &  0.011 &    21 &  0.582 &     31 &  0.355 \\
     9 &       4275 & 204748\_at & PTGS2  &  18 &  14 &   424 &  0.118 &    13 &  0.624 &     42 &  0.313 \\
    10 &       1437 & 201909\_at & RPS4Y1  &  21 &  34 &  8232 &  0.013 &     3 &  0.913 &      5 &  0.736 \\
    11 &      14150 & 214774\_x\_at & TOX3  &  34 &   0 &    95 &  0.262 &    49 &  0.492 &    102 &  0.210 \\
    12 &      10518 & 211074\_at & FOLR1  &   7 &   4 &  9482 &  0.011 &   926 &  0.124 &    213 &  0.159 \\
    13 &       9580 & 210096\_at & CYP4B1  &   6 &  44 &     8 &  0.797 &    65 &  0.431 &     54 &  0.284 \\
    14 &       4002 & 204475\_at & MMP1 &  27 &   0 &    34 &  0.406 &    24 &  0.564 &     21 &  0.465 \\
    15 &      13990 & 214612\_x\_at & MAGEA  &  16 &   0 &   489 &  0.110 &   134 &  0.323 &     35 &  0.339 
    \end{tabular}\end{center}
    \vspace*{-18pt}
\end{table}

\begin{table}[t!]
\caption{\label{tbl:gene_ls2}
Genetics example: the probes with top (scaled) leverage scores, derived from the first two singular vectors.}

\begin{center}\small
\begin{tabular}{crr|ll|cc|c} \hline
\multicolumn{1}{c}{LS (2) } &
\multicolumn{1}{c}{index} &
\multicolumn{1}{c}{LS (2)} &
\multicolumn{1}{|c}{probe} &
\multicolumn{1}{c}{gene} &
\multicolumn{1}{|c}{number} &
\multicolumn{1}{c}{number} &
\multicolumn{1}{|c}{DEIM } \\
\multicolumn{1}{c}{rank } &
\multicolumn{1}{c}{$q_j$} &
\multicolumn{1}{c}{score} &
\multicolumn{1}{|c}{set} &
\multicolumn{1}{c}{name} &
\multicolumn{1}{|c}{sick } &
\multicolumn{1}{c}{well} &
\multicolumn{1}{|c}{rank } \\ \hline
     1 &  9565 &  1.000 & 210081\_at & AGER &   2 &  45 &    1  \\ 
     2 & 13766 &  0.922 & 214387\_x\_at & SFTPC&   6 &  48 &  ---  \\ 
     3 & 11135 &  0.907 & 211735\_x\_at & SFTPC&   5 &  48 &  73  \\ 
     4 &  9361 &  0.899 & 209875\_s\_at & SPP1&  50 &   2 &  ---  \\ 
     5 &  5509 &  0.896 & 205982\_x\_at & SFTPC&   5 &  48 &  ---  \\ 
     6 &  9103 &  0.835 & 209613\_s\_at & ADH1B&   2 &  47 &  ---  \\ 
     7 & 14827 &  0.834 & 215454\_x\_at & SFTPC&   0 &  46 &  ---  \\ 
     8 &  9580 &  0.797 & 210096\_at & CYP4B1&   6 &  44 &   13  \\ 
     9 &  4239 &  0.754 & 204712\_at & WIF1&   5 &  43 &  70  \\ 
    10 &  3507 &  0.724 & 203980\_at & FABP4&   2 &  44 &  ---  \\ 
    11 & 18594 &  0.717 & 219230\_at & TMEM100&   2 &  38 &  ---  \\ 
    12 &  9102 &  0.684 & 209612\_s\_at & ADH1B&   2 &  46 &  ---  \\ 
    13 & 13514 &  0.626 & 214135\_at & CLDN18&   3 &  47 &  ---  \\ 
    14 &  5393 &  0.626 & 205866\_at & FCN3&   0 &  39 &  ---  \\ 
    15 &  4727 &  0.614 & 205200\_at & CLEC3B&   0 &  39 &  ---  \\ 
\end{tabular}\end{center}
    \vspace*{-18pt}

    \end{table}

While the rows selected from leverage scores did not produce as accurate an approximation, $\|\bA-\bC_k\bU_k\bR_k\|$, as DEIM, these probes do a much more effective job of discriminating patients with tumors from those without.  Indeed, for~14 of the top~15 probes, the tumor-free patients express strongly, while the patients with tumors do not;
in the remaining case, the opposite occurs.

\vspace*{-10pt}

\section{Conclusions}  

\vspace*{-5pt}

\noindent
The Discrete Empirical Interpolation Method (DEIM) is an index selection procedure that gives simple, deterministic CUR factorizations of the matrix $\bA$.\ \ 
Since DEIM utilizes (approximate) singular vectors, we propose an effective one-pass incremental approximate QR factorization that can efficiently compute dominant singular vectors for data sets with rapidly decaying singular values; this method could prove useful in a variety of other settings.    
The accuracy of the resulting rank-$k$ CUR factorization can be bounded in terms of $\sigma_{k+1}$, the error in the best rank-$k$ approximation to $\bA$.\ \  Our analysis of the 2-norm error $\|\bA-\bC\bU\bR\|$ applies to all CUR approximations that use the optimal central factor $\bU=\bC^I\bA\bR^I$, and hence can give insight into the performance of other index selection algorithms, such as leverage scores, uniform random sampling, or column-pivoted QR factorization.  Numerical examples illustrate that the DEIM-CUR approach can deliver very good low rank approximations, compared to row selection based on dominant leverage scores.

\vspace*{-7pt}
\section*{Acknowledgements}

\vspace*{-7pt}
We thank Inderjit Dhillon, Petros Drineas, Ilse Ipsen, Michael Mahoney, 
and Nick Trefethen for a number of helpful discussions.  
We are also grateful to Gunnar Martinsson for recommending experiments with the column-pivoted QR selection algorithm, and to an anonymous referee for encouraging us to seek the improved analysis and growth example for the DEIM error constants at the end of Section~\ref{sec:theory}.



\bibliographystyle{plain}
\bibliography{CUR.bib}

\end{document}